\documentclass[10pt]{iopart}
\usepackage{iopams}
\usepackage{epsfig}

\usepackage{epsfig}

\usepackage{xcolor}
\usepackage[active]{srcltx}

\begin{document}
\newtheorem{col}{Corollary}
\newtheorem{lm}{Lemma}
\newtheorem{theorem}{Theorem}
\newtheorem{df}{Definition}
\newtheorem{prop}{Proposition}
\newtheorem{remark}{Remark}
\eqnobysec

\def\fl{\!}
\def\TT{{\cal T}}
\def\ds{\displaystyle}
\def\proof{\noindent {\sf Proof.}  }
\def\qed{\hfill $\Box$ \\ \bigskip}
\def\Id{\mbox{Id}}
\def\loc{\mbox{loc}}

\newcommand{\N}{\ensuremath{\mathbb{N}}}
\newcommand{\Z}{\ensuremath{\mathbb{Z}}}
\newcommand{\Q}{\ensuremath{\mathbb{Q}}}
\newcommand{\R}{\ensuremath{\mathbb{R}}}
\newcommand{\C}{\ensuremath{\mathbb{C}}}
\newcommand{\SI}{\ensuremath{\mathbb{S}^1}}
\newcommand{\T}{\ensuremath{\mathbb{T}}}


\newcommand{\Rrev}{R}   
\newcommand{\Lrev}{L}   
\newcommand{\Wu}{W^u}  
\newcommand{\Ws}{W^s}  
\newcommand{\sech}{\mathrm{sech}}

\newcommand{\Mt}{\tilde{M}}
\newcommand{\ct}{\tilde{c}}

\hyphenation{dif-feo-mor-phisms}


\title[Attracting, repelling and elliptic orbits in reversible maps]
{Abundance of attracting, repelling and elliptic periodic orbits
in  two-dimensional reversible maps.}

%
%

\author{Delshams A.\dag\ , Gonchenko S.V.\ddag , Gonchenko
V.S.\ddag ,
L\'azaro J.T.\dag\ and Sten'kin O.\ddag }

\address{\dag\ Universitat Polit\`ecnica de Catalunya, Barcelona, Spain}

\address{\ddag\ Institute for Applied Mathematics \& Cybernetics, N.Novgorod, Russia}

\begin{abstract}
We study dynamics and bifurcations of two-dimensional reversible
maps having non-transversal heteroclinic cycles containing
symmetric saddle periodic points. We consider one-parameter
families of reversible maps unfolding generally the initial
heteroclinic tangency and prove that there are infinitely
sequences (cascades) of bifurcations of birth of asymptotically
stable and unstable as well as elliptic periodic orbits.

\end{abstract}


\section{Introduction}
Reversible systems have a very special status inside the realm of dynamical
systems. Usually, they are positioned ``between'' dissipative and conservative
systems. In the context of continuous dynamical systems, a reversibility means that the system is invariant under
the change of time-direction, $t\mapsto -t$, and a transformation in the spatial variables. In the discrete context, reversibility of a map $f$ (a diffeomorphism) means
that $f$ and $f^{-1}$ possess the same dynamics. Notice that
the term ``the same'' can have rather different meaning. If $f$ and $f^{-1}$
are smoothly conjugate, i.e., $f\circ h = h\circ f^{-1}$ and $h$ is a just a
diffeomorphism, then $f$ is called \emph{weekly reversible}. However, much more
interesting types of reversibility appear when $h$ possesses some structures.
For example, if $h$ is an \emph{involution}, i.e., $h^2=\mathrm{Id}$. In this
case, the map $f$ is called \emph{strongly reversible}. Since this last case is
the most frequent one in the literature (probably, beginning with Birkhoff),
nowadays strongly reversible maps are simply called \emph{reversible} maps.

In contrast to conservative and dissipative systems, the study of homoclinic
bifurcations in reversible systems is not so popular. Even for two-dimensional
maps, only few results are known and most of them relate to ``conservative and
reversible'' maps which form a certain codimension-$\infty$ subclass in the
class of reversible maps. This situation is probably due to the ``common
belief'' that conservative and dissipative phenomena of dynamics only exist
separately and, thus, there is no necessity to study them ``all together''.

However, they actually can appear together in a dynamical system, giving rise
to the so-called \emph{phenomenon of mixed dynamics}, which was recently
discovered in \cite{GST97} (see also \cite{GStS02,GStS06,VGSh07}). The essence
of this phenomenon consists in the fact that
\begin{enumerate}
\item[(i)] a dynamical system has simultaneously infinitely many hyperbolic
    periodic orbits of all possible types (stable, completely unstable and
    saddle), and
\item[(ii)] these orbits are not separated as a whole, i.e., the closures
    of sets of orbits of different types have nonempty intersections.
\end{enumerate}

It was shown in \cite{GST97} that the property of mixed dynamics can be
generic, i.e., it holds for residual subsets of open regions of systems. In
particular, it was also proved that such regions (Newhouse regions,
in fact) exist near two-dimensional diffeomorphisms with non-transversal
heteroclinic cycles containing at least two saddle periodic points $O_1$ and
$O_2$ such that $|J(O_1)|>1$ and $|J(O_2)|<1$, where $J(O_i)$ is the Jacobian
of the Poincar\'{e} map (the diffeomorphism iterated as many times as the period of
$O_i$) at the point $O_i$, $i=1,2$.

Let us recall that a \emph{heteroclinic cycle} (contour) is a set consisting of
saddle hyperbolic periodic orbits $O_1,\dots,O_n$ as well as heteroclinic
orbits $\Gamma_{i,j}\subset W^u(O_i)\cap W^s(O_j))$, where at least the orbits
$\Gamma_{i,i+1}$ and $\Gamma_{n,1}$ for $i=1,\dots,n-1$, are included. In
general, cycles can include also homoclinic orbits $\Gamma_{i,i}\subset
W^u(O_i)\cap W^s(O_i)$. An heteroclinic cycle is called \emph{non-transversal}
(or non-rough) if at least one of the pointed out intersections $W^u(O_i)\cap
W^s(O_j)$ is not transverse.

If a heteroclinic cycle (or a homoclinic orbit) is transverse, then, as is
well-known after Shilnikov \cite{Sh67}, the set of orbits entirely lying in a
small neighbourhood is a locally maximal uniformly hyperbolic set. The situation
becomes \emph{drastically} different in the non-transversal case. One can say
even that the corresponding system is infinitely degenerate, since its
bifurcations can produce homoclinic tangencies of arbitrary high orders and, as
a consequence, arbitrary degenerate periodic orbits \cite{GST93a,GST99,GST07}.

We remind also that systems with homoclinic tangencies are dense in open
regions (the so-called \emph{Newhouse regions}) in the space of smooth
dynamical systems \cite{N69,N74,N79}. Moreover, these regions exist near any
system with a homoclinic tangency (or a non-transversal heteroclinic cycle).
Importantly, these regions are present in parameter families unfolding
generally the initial homoclinic (or heteroclinic) tangency in certain open
domains of the parameter space in which there are dense values of the
parameters corresponding to the existence of homoclinic tangencies. Certainly,
such domains are called again Newhouse (parameter) regions or  Newhouse
intervals for one-parameter families. In general, it should be clear from the
context the kind of Newhouse regions considered.

The existence of Newhouse regions near systems with homoclinic tangencies was
established in \cite{N79} for two-dimensional diffeomorphisms, in
\cite{GST93b,PV94,Romero95} for the general multidimensional case (including
parameter families \cite{N79,GST93b}) and in \cite{D00} for area-preserving
maps. The existence of Newhouse regions near systems with non-transversal
heteroclinic cycles follows immediately from these results, since in such case
homoclinic tangencies appear under arbitrary small perturbations. Moreover, in
this case also the so-called \emph{Newhouse regions with heteroclinic
tangencies} can exist. It was proved in \cite{GST97} that if the non-transversal
heteroclinic cycle is simple, i.e., it contains only one non-transversal
heteroclinic orbit and the corresponding tangency of invariant manifolds is
quadratic, then Newhouse intervals with heteroclinic tangencies exist in any
general one-parameter unfolding.

\begin{figure}[bht]
\begin{center}
\includegraphics[width=12cm]{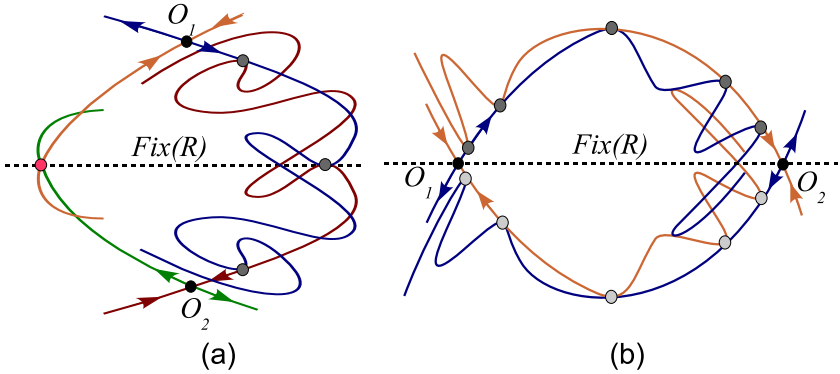}
\caption{{\footnotesize Two examples of planar reversible maps
with symmetric non-transversal  (quadratic tangency) heteroclinic
cycles. One can think, for simplicity, that the involution $R$ in
both cases is linear of form  $x\mapsto x, y\mapsto -y$ and, thus,
$\mbox{Fix}\;R = \{y=0\}$.
}} \label{fig:Intro1}
\end{center}
\end{figure}

The above-mentioned mixed dynamics takes place as a generic phenomenon. This is the case, for instance, when the initial heteroclinic cycle is \emph{contracting-expanding}~\cite{GST97},
that is, when it contains contracting and expanding periodic points (i.e. with
the absolute value of its Jacobian being greater or less than $1$). It is
worth mentioning that contracting-expanding heteroclinic cycles are rather
usual among reversible maps. An example of such a cycle is shown in
Figure~\ref{fig:Intro1}(a). In this example the reversible map has two saddle
fixed points $O_1$ and $O_2$ and two heteroclinic orbits $\Gamma_{12}\subset
W^u(O_1)\cap W^s(O_2)$ and $\Gamma_{21}\subset W^u(O_2)\cap W^s(O_1)$ such that
$R(O_1)=O_2$ and $R(\Gamma_{21})=\Gamma_{21}$, $R(\Gamma_{12})=\Gamma_{12}$.
Besides, the orbit $\Gamma_{12}$ is non-transversal, so that the manifolds
$W^u(O_1)$ and $W^s(O_2)$ have a quadratic tangency along $\Gamma_{12}$.
Since $R(O_1)=O_2$, it turns out that their Jacobians verify $J(O_1)=J^{-1}(O_2)$. If $J(O_i) \neq \pm 1$, $i=1,2$, then the heteroclinic
cycle is contracting-expanding.  This condition is robust and is
perfectly compatible with reversibility.

Certainly, results of \cite{GST97} can be applied to reversible maps with such
heteroclinic cycles and so the phenomenon of mixed dynamics becomes very important
and generic. However, reversible systems are sharply different from general
ones by the fact that they can possess robust non-hyperbolic symmetric periodic
orbits, more precisely, \emph{elliptic symmetric periodic points}. Thus, one
realizes that the phenomenon of mixed dynamics in the case of two-dimensional
reversible maps  should be connected with the
\emph{coexistence of infinitely many attracting, repelling, saddle and elliptic
periodic orbits}.
The existence of Newhouse regions (intervals) in which this property is generic
was already established in \cite{LSt} for the case of reversible
two-dimensional maps close to a map having a heteroclinic cycle of the type
depicted in Figure~\ref{fig:Intro1}(a).

However, it appears to be true that the phenomenon of mixed dynamics is
universal for reversible (two-dimensional) maps with complicated dynamics when
symmetric structures (symmetric periodic, homoclinic and heteroclinic orbits)
are involved. This universality can be formulated as the following

 \smallskip
\noindent{\bf Reversible Mixed Dynamics Conjecture\ }
\emph{ Two-dimensional reversible maps with mixed dynamics are generic (compose residual subsets) in Newhouse regions in which there are dense maps with symmetric homoclinic or/and heteroclinic tangencies}.

 \smallskip

We will assume, in what follows, that the involution $R$
is not trivial, i.e. it satisfies
\begin{equation}
\label{invreq}
R^2 = \mbox{Id},\;\; \mbox{dim}\;\mbox{Fix}\;R = 1.
\end{equation}

We will say that an object $\Lambda$ is \emph{symmetric} when $R(\Lambda)=\Lambda$. To put more
emphasis, sometimes the notation self-symmetric may be used. By a
\emph{symmetric couple of objects} $\Lambda_1,\Lambda_2$, we will mean two different
objects that are symmetric to each other, i.e., $R(\Lambda_1)=\Lambda_2$.

Then the \emph{symmetric homoclinic (heteroclinic) tangencies} from the
RMD-Conjecture can be divided into two main types: 1) there is a
non-transversal symmetric heteroclinic orbit to a symmetric couple of saddle
points, or 2) there is a symmetric couple of non-transversal homo/heteroclinic
orbits to symmetric saddle points.

Notice that the heteroclinic quadratic tangency shown in
Fig.~\ref{fig:Intro1}(a) relates to the type~1), whereas an example of
reversible map having a heteroclinic cycle of type~2) is shown in
Figure~\ref{fig:Intro1}(b). The latter map has two symmetric saddle fixed
points $O_1$ and $O_2$ ($R(O_1)=O_1$, $R(O_2)=O_2$) and a symmetric couple of
heteroclinic orbits $\Gamma_{12}\subset W^u(O_1)\cap W^s(O_2)$ and
$\Gamma_{21}\subset W^u(O_2)\cap W^s(O_1)$ ($R(\Gamma_{12})=\Gamma_{21}$).

As we mentioned above, the case of non-transversal heteroclinic cycles of
type~1), as in Fig.~\ref{fig:Intro1}(a), was studied in the paper \cite{LSt}
where, in fact, the  RMD-Conjecture was proved for general one-parameter
(reversible) unfoldings, under the generic condition $J(O_1)=J^{-1}(O_2) \neq
1$.

In the first part of this paper, Sections~\ref{se:2} and~\ref{sec:prTh1}, we
state and prove the RMD-Conjecture for one-parameter families which unfold
generally heteroclinic tangencies of type~2), as in Fig.~\ref{fig:Intro1}(b).
We will call this type as \emph{reversible maps with a symmetric couple of
heteroclinic tangencies}. We notice that in this case the condition
$J(O_1)=J(O_2)=1$ holds always since $O_i\in \mbox{Fix}\; R$, showing that
generic conditions are different in systems with different types. The generic
condition that we are going to assume for systems of type~2) is denoted by
condition~\textsf{[C]} in Theorem~\ref{th:main1}. This condition
(see~(\ref{condition:C}) and comments to it) amounts to say that the (global)
map defined near a heteroclinic point is neither a uniform contraction
(expansion) nor a conservative map. More precisely, it has to have a
non-constant Jacobian in those local coordinates (near $O_1$ and $O_1$) in
which the saddle maps are, a priori,  area-preserving. In particular, such local
coordinates are given by Lemma~\ref{LmLocalMap} in which the normal form of the
first order for a saddle map is derived.\footnote{However, the property of a
symmetric saddle periodic point to be a priori area-preserving is more
delicate. It is well-known, see e.g. \cite{DL05}, that a symmetric reversible
saddle map is ``almost conservative'', i.e. its analytical normal form is
exactly conservative, and its $C^\infty$ formal normal form (up to ``flat
terms") is conservative.)}

Moreover, as we show, symmetry breaking bifurcations have also another nature,
in comparison with \cite{LSt}. We find a two-step
``fold$\Rightarrow$pitch-fork'' scenario of bifurcations in the first-return
maps leading to the appearance of non-conservative fixed points which can be
either attracting and repelling or saddle with Jacobian greater and less than
$1$ (see Theorem~\ref{th:main1} and Figure~\ref{bifseq}). Notice that in the
case of heteroclinic cycles of type~1), such non-conservative points appear
just under fold bifurcations \cite{LSt} in a symmetric couple of first-return
maps, whereas elliptic points appear in symmetric first-return maps.

\begin{figure}[bht]
\begin{center}
\includegraphics[height=8cm]{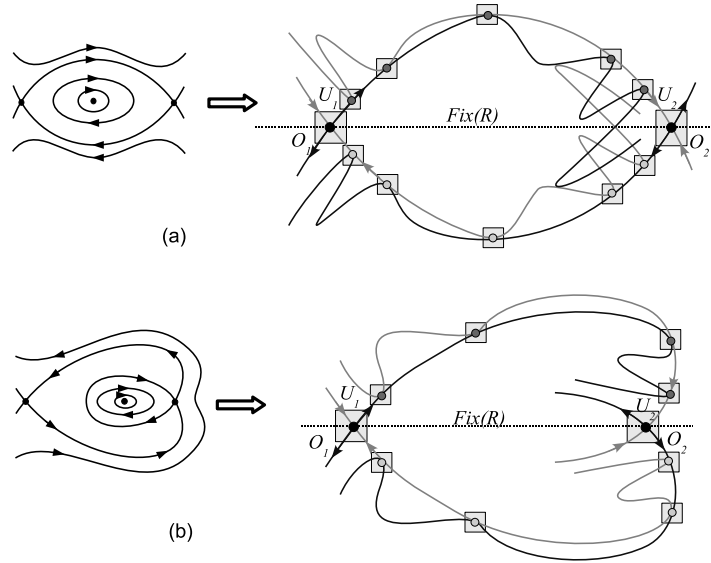}
\caption{Examples of reversible maps with non-transversal heteroclinic cycles in
the cases of  (a) ``inner tangency''; (b) ``outer tangency''. These diffeomorphisms can be
constructed as the Poincar\'e maps for periodically perturbed conservative planar systems (left) when
the corresponding perturbations preserve the reversibility (for example, with respect to the change $x\to x,y\to -y,t\to -t$,
as in the system $\dot x = y$, $\dot y = x - x^3 + \varepsilon \dot x \cos t$). Notice that the resulting reversible map cannot be area-preserving, in general}
\label{fig:2class}
\end{center}
\end{figure}

The second part of this work, Section~\ref{se:examples}, has a more applied
character. We show in Subsection~\ref{subse:Duffing} that reversible
two-dimensional maps with a priori non-conservative orbit behaviour can
be obtained as certain periodic perturbations of two-dimensional conservative
flows of form $\dot x =y,\; \dot y = F(y)$ . We require that these
perturbations include explicitly the ``friction term'' $\dot x$ and preserve
only reversible properties of the initial flow (for example, they keep the
perturbed systems to be invariant under the change $x\mapsto x, y\mapsto -y,
t\mapsto -t$). In this way, in particular, we can obtain the reversible maps of
type~2), see Figure~\ref{fig:2class}. As a concrete example, we consider, in
Subsection~\ref{subse:Duffing}, the periodically perturbed Duffing equation. In
Subsection~\ref{subse:PikTopaj} we consider an example of a reversible map from
\cite{PT02} defined on the two-dimensional torus which shows a visible
non-conservative orbit behaviour. This map describes the dynamics of three
coupled simple rotators with small symmetric couplings. The couplings are
chosen here in such a way that they preserve the reversibility of the initial
uncoupled three simple rotators.

In the third part of the paper, at Section~\ref{se:crossL1L2}, we consider a series of problems related to the
representation of reversible two-dimensional maps in the so-called {\em
cross-form}. It is well-known that such type of cross-forms of maps are very
convenient for studying hyperbolic properties of systems with homoclinic
orbits, both transversal~\cite{Sh67} and
non-transversal~\cite{GaS73,G84,GST08}.\footnote{Since L.P. Shilnikov was the
first author introducing such forms and coordinates into dynamical systems,
they are often referred to as  ``Shilnikov cross-form'' and ``Shilnikov
coordinates''.} We show that the cross-forms and the corresponding
cross-coordinates are natural for reversible maps, since they allow expressing
many reversible structures explicitly and simplify analytical treatment. In
particular, in the corresponding local cross-coordinates near a symmetric
saddle periodic point, the normal form of the saddle map becomes very simple.
And last (but not least), Section~\ref{se:technical_lemmas}
contains the proofs of the lemmas needed in the proof of the main theorem~\ref{th:main1}.

\subsection{Out of the general rule: a collection of reversible maps with codimension one homoclinic and heteroclinic tangencies}
\begin{figure}[bht]
\begin{center}
\includegraphics[width=12cm]{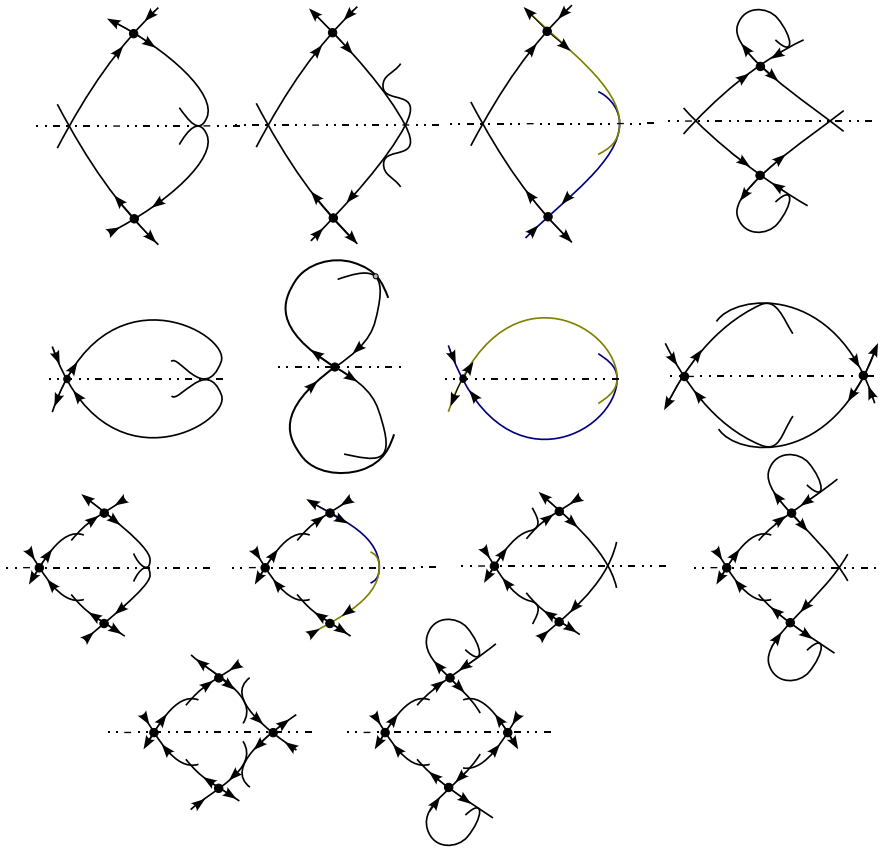}
\caption{{\footnotesize Some examples of reversible maps with homoclinic and heteroclinic tangencies}} \label{fig:Intro1-typs}
\end{center}
\end{figure}

It is important to notice that there are  many other cases of reversible maps
with homoclinic and heteroclinic tangencies for which one needs to prove the
RMD-Conjecture in the framework of one-parameter general families. In
Figure~\ref{fig:Intro1-typs} we collect some simple examples of such maps. They
differ by the type of fixed points and tangencies: homoclinic or heteroclinic,
quadratic or cubic\footnote{The existence of symmetric cubic homoclinic or
heteroclinic tangencies is a codimension one bifurcation phenomenon in the
class of reversible maps. Therefore, these cubic tangencies should be also
considered as the main ones jointly with the pointed out quadratic ones.}, etc.
However, it turns out to be more important that all this kind of maps could be
separated into two groups: a first one including those maps which are a priori non-conservative  and a second one with those maps where this
non-conservativity is, in some sense, hidden.

A map as in Figure~\ref{fig:Intro1}(a) (the first one in Figure~\ref{fig:Intro1-typs}) belongs to the first group, since the
condition $J(O_{1,2})\neq 1$ destroys certainly the conservative character.
Indeed, under splitting such heteroclinic cycle, homoclinic tangencies appear
both  to saddles with Jacobian greater and less than 1 and, thus, attracting
and repelling periodic points can be born \cite{GaS73}. By this principle, all
the other cases where a symmetric couple of fixed (periodic) points are
involved, can be referred to as the class of a priori non-conservative maps
(e.g., all maps of the first, third and fourth rows of
Figure~\ref{fig:Intro1-typs}). For maps of this type, the problem of finding
symmetric periodic orbits (i.e., elliptic ones) has to be considered as very
important.

The maps at the second row of
Figure~\ref{fig:Intro1-typs} have only symmetric fixed points.
Evidently, the first, third and fourth maps can be assigned to maps
with "hidden non-conservativity", since it is not clear, in advance, the existence of
bifurcation mechanisms leading to the appearance of attracting and
repelling periodic orbits. It is not the case of the third map (at the
second row) which has a
symmetric couple of (quadratic) homoclinic tangencies to the same symmetric fixed
saddle point. One can assume, without loss of reversibility,
that the map near a homoclinic point is not
conservative (the Jacobian is greater or less than $1$). Then,
clearly,  stable (unstable) periodic orbits can be born under such
homoclinic bifurcations. Thus, the main problem here is to prove the
appearance of elliptic periodic orbits and, as a first step, to do it in the
one-parameter setting.

\subsection{A short description of the main results}
\label{subse:Main_Theorems}

We describe briefly the central ideas underlying our main results (Theorems~\ref{th:main1} and~\ref{th:col-new}). Let assume $f_0$ to be a $R$-reversible map of type as in Figure~\ref{fig:Intro1}(b), that is, having two symmetric
saddle fixed points $O_1$ and $O_2$  (i.e. $R(O_i)=O_i$, $i=1,2$)
and two asymmetric non-transversal heteroclinic orbits
$\Gamma_{12}\subset W^u(O_1)\cap W^s(O_2)$ and $\Gamma_{21}\subset W^u(O_2)\cap W^s(O_1)$, satisfying that $R(\Gamma_{12})=\Gamma_{21}$. Let us consider $\left\{ f_{\mu} \right\}$ any general one-parameter unfolding of $f_0$ with $\mu$ being the parameter splitting the initial heteroclinic tangency. The main goal of this paper is to show that under general hypotheses, any of these unfoldings undergoes infinitely many (in fact, a cascade) of symmetry-breaking bifurcations of single-round periodic orbits.
Such bifurcations in these first-return maps (defined near some point of
the heteroclinic tangency) follow from the following scenario:
\footnote{Notice that in the Lamb-Sten'kin case~\cite{LSt}, see
Fig.~\ref{fig:Intro1}(a), fold bifurcations are directly symmetry breaking and lead to the
appearance of two pairs of asymmetric periodic orbits: (saddle, sink) and (saddle, source).}
\[
 \mbox{fold\;bifurcation}\;\;\Rightarrow\;\;  \mbox{pitch-fork\;
bifurcation}
\]
Under fold-bifurcation, a symmetric parabolic fixed point appears which falls afterwards into two symmetric saddle and elliptic fixed points. Concerning (reversible) pitch-fork bifurcations, they can be of two classes depending essentially on the type of initial heteroclinic cycle they exhibit (see Figure~\ref{fig:2class} and \ref{bifseq}):
in case (a) we say that $f_0$ has a heteroclinic cycle of ``inner tangency'' while in case (b) we say that it is of ``outer tangency''.
In the first case, ``inner tangency'',  from the
symmetric elliptic fixed point three fixed points are born: a symmetric saddle and asymmetric sink and source.
In the second case, ``outer tangency'',
under a pitch-fork bifurcation, the symmetric saddle point becomes a symmetric elliptic
point and two fixed asymmetric saddle points with Jacobian greater and less than 1.
Both scenarios are showed in Figure~\ref{bifseq}.
In the bifurcation diagram of Figure~\ref{Figdiag31},
related to a conservative approximation of the rescaled first-return map,  value of $\tilde c<0$ corresponds to the ``inner tangency'' case and $\tilde c>0$ corresponds to the ``outer tangency'' one
(situation $\tilde c=0$ is singular and not realisable for our maps; therefore there are no transitions from $\tilde c<0$ to $\tilde c>0$).

We finish this Introduction by describing the structure of the paper.
Section~\ref{se:2} is devoted to the two above-mentioned main results of the
paper, Theorems~\ref{th:main1} and~\ref{th:col-new}. Their proof is presented
in Section~\ref{sec:prTh1} and relies to five lemmas whose proof is deferred to
Sections~\ref{se:crossL1L2}, \ref{se:technical_lemmas} and~\ref{se:birkhoff}.
Section~\ref{se:examples} contains two concrete examples of applications of the
main results of this paper.

\section{Symmetry breaking bifurcations in the case of reversible maps with non-transversal heteroclinic cycles}
\label{se:2}

Let $f_0$ be a $C^r$-smooth, $r\geq4$, two-dimensional map, reversible with respect to an involution $R$ satisfying $\dim\;\mbox{Fix}(R)=1$. Let us assume that $f_0$ satisfies the following two conditions:
\begin{itemize}
\item[\textsf{[A]}]
$f_0$ has two saddle fixed points $O_1$ and $O_2$ belonging to the line $\mbox{Fix}(R)$
and that any point $O_i$ has multipliers $\lambda_i,\lambda_i^{-1}$ with
$0<\lambda_i<1$, $i=1,2$.

\item[\textsf{[B]}]
The invariant manifolds $W^u(O_1)$ and $W^s(O_2)$ have quadratic tangencies at the points of some heteroclinic orbit $\Gamma_{12}$ and, therefore, by reversibility, the manifolds $W^u(O_2)$ and $W^s(O_1)$ have quadratic tangencies at the points of a heteroclinic orbit $\Gamma_{21} =R(\Gamma_{12})$.
\end{itemize}
Hypotheses \textsf{[A]-[B]}  define reversible maps with non-transversal symmetric heteroclinic cycles like in Figure~\ref{fig:Intro1}(b).
We ask them to satisfy one more condition.
Namely, consider
two points $M_1\in W^u_{loc}(O_1)$ and $M_2\in W^s_{loc}(O_2)$ belonging to the same heteroclinic orbit $\Gamma_{12}$ and suppose $f_0^q(M_1)=M_2$ for a suitable integer $q$. Let some smooth local coordinates
$(x_i,y_i)$ be chosen near the points $O_i$ in such a way that the local invariant manifolds are straightened, i.e., $W^u_{loc}(O_i)$ and $W^s_{loc}(O_i)$ have, respectively, equations $x_i=0$ and $y_i=0$. Let $T_{12}$ denote the
restriction of the map $f_0^q$ onto a small neighbourhood of the point $M_1$. Then, we assume that
\begin{itemize}
\item[\textsf{[C]}]
the Jacobian of $T_{12}$ is not constant and, moreover,
\begin{equation}
\label{condition:C}
Q= \frac{\partial J(T_{12})}{\partial y}\bigl|_{M_1}\neq 0
\end{equation}
\end{itemize}

Condition $J(T_{12})\neq \;\mbox{const}$ is well defined only when certain
restrictions on the local coordinates hold. One possibility
is when these coordinates $(x_i,y_i)$ around $O_i$ are chosen in such a way that $W^u_{loc}(O_1)$ and $W^s_{loc}(O_2)$ are straightened. However, the sign of $J(T_{12})$ depends also on the orientation chosen for the coordinate axes.
To be precise, we choose these orientations in such a a way that: $(i)$ the $y$ and
$x-$coordinates of the heteroclinic points $M_1\in W^u_{loc}(O_1)$ and $M_2\in W^s_{loc}(O_2)$ are positive; $(ii)$ for the symmetric points, $M_1^\prime = R(M_1)\in W^s_{loc}(O_1)$ and $M_2^\prime = R(M_2)\in W^u_{loc}(O_2)$, the $x$ and $y$-coordinates are positive as well.

Two classes of reversible maps satisfying conditions \textsf{[A]-[C]}
can be distinguished: those maps with {\em
``inner''} (heteroclinic) tangency and those with {\em ``outer''}  tangency, corresponding  to $J(T_{12})>0$ and $J(T_{12})<0$, respectively.
Two examples of such diffeomorphisms are shown in Figure~\ref{fig:2class}. Notice that in both cases the global map $T_{12}$ is orientable.
In the case (a) the axes $x_1,y_1$ and $x_2,y_2$ have the same orientation, whereas the orientations are different in the case (b).

Once stated the general conditions for $f_0$, let us
embed it into a one-parameter family $\left\{ f_{\mu} \right\}$ of reversible maps that
unfolds generally at $\mu=0$ the initial heteroclinic
tangencies at the points of $\Gamma_{12}$. Then, without loss of
generality, we can take $\mu$ as the corresponding splitting parameter.
By reversibility, the invariant
manifolds $W^u(O_1)$ and $W^s(O_2)$ split as $W^u(O_2)$ and
$W^s(O_1)$ do when $\mu$ varies. Therefore, since these
heteroclinic tangencies are quadratic, only one governing parameter
is needed to control this splitting.

\begin{figure}[bht]
\begin{center}
\includegraphics[width=12cm]{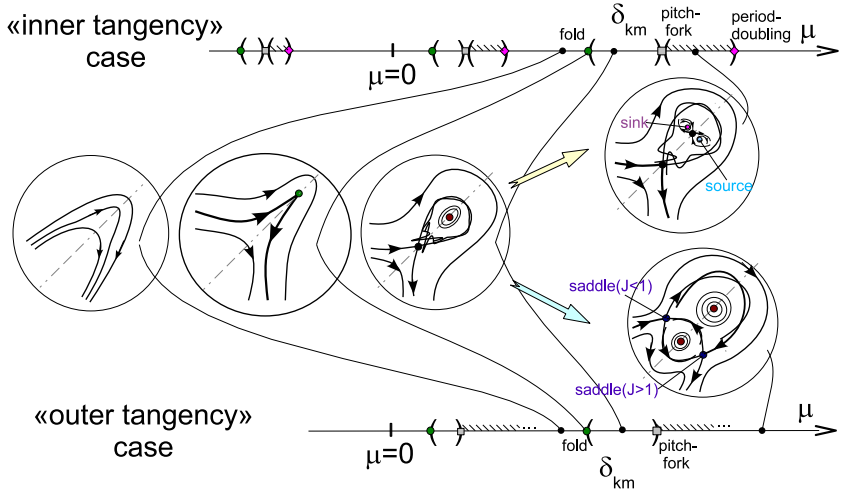} \caption{We mark, by shading, some intervals corresponding to the existence of two
asymmetric single-round periodic orbits.}
\label{bifseq}
\end{center}
\end{figure}

Let $U$ be an small enough neighbourhood of the contour
$C=\{O_1,O_2,\Gamma_{12},\Gamma_{21}\}$. It can be represented as the union of two small
neighbourhoods (disks) $U_1$ and $U_2$ of the saddles $O_1$ and $O_2$ and a finite number of small disks containing those points of $\Gamma_{12}$ and $\Gamma_{21}$ which do not belong to $U_1$ and $U_2$ (see Figure~\ref{fig:2class}).
We will focus our attention on the bifurcations of the so-called {\em single-round periodic orbits}, that is, orbits lying entirely in $U$ and
having exactly one intersection point with every disk from the set $U\backslash(U_1\cup U_2)$. Any point of a single-round periodic orbit is a fixed point of the corresponding {\em first-return map} $T_{km}$, that is constructed by orbits of $f_\mu$ with $k$ and $m$  iterations (of $f_\mu$) in $U_1$ and $U_2$, respectively.
We will call them single-round periodic orbit of type $(k,m)$. The values of $k$ and $m$ will be always prescribed \emph{a priori}. The first main result is as follows:
\begin{theorem}
Let $\left\{f_\mu\right\}_{\mu}$ be a one-parameter family of reversible diffeomorphisms that unfolds, generally, at $\mu=0$ the initial heteroclinic tangencies. Assume that $f_0$ verifies the conditions {\rm \textsf{[A]-[C]}}.

Then, at any segment
$[-\epsilon,\epsilon]$ with $\epsilon>0$ , there are infinitely many intervals $\delta_{km}$ with border points $\mu_{\rm fold}^{(k,m)}$ and $\mu_{\rm pf}^{(k,m)}$
such that $\delta_{km}\to 0$ as $k,m\to\infty$ and the following
holds:
\begin{itemize}

\item[$(i)$] The value $\mu=\mu_{\rm fold}^{(k,m)}$ corresponds to a
non-degenerate conservative fold bifurcation and, thus, the
diffeomorphism $f_\mu$ has at $\mu\in\delta_{km}$ two symmetric,
saddle and elliptic, single-round periodic orbits of type $(k,m)$.

\item[$(ii)$] The value $\mu = \mu_{\rm pf}^{(k,m)}$ corresponds to a symmetric (and non-degenerate if condition {\sf [C]} holds) pitch-fork bifurcation depending on the type of
$f_0$:

\begin{itemize}
\item[$(ii)_a$]  In the case of ``inner'' tangency, single-round asymmetric attracting and repelling periodic orbits of type $(k,m)$ are born and, moreover, these orbits undergo simultaneously non-degenerate period doubling bifurcations at the value
$\mu = \mu_{\rm pd}^{(k,m)}$ (where $\mu_{\rm pd}^{(k,m)}\to 0$ as $k,m\to\infty$).

\item[$(ii)_b$] For the ``outer'' tangency,
two single-round saddle periodic orbits of type $(k,m)$ with Jacobian greater and less than $1$, respectively, are born. Moreover, they do not bifurcate any more (at least for $|\mu|<\epsilon$).
\end{itemize}
\end{itemize}
\label{th:main1}
\end{theorem}
We refer the reader to Figure~\ref{bifseq} for an illustration of this theorem.

Theorem~\ref{th:main1} and its counterpart result in \cite{LSt} show that the appearance of non-conservative periodic orbits under global bifurcations can be consider as a certain generic property of two-dimensional reversible maps.

Briefly, the method we use - based on a rescaling technique - will allow us to prove
that the first-return map $T_{km}$ can be written asymptotically close (as $k,m\to\infty$)
to an area-preserving map of the form:
\begin{equation}
H\;:\left\{
\begin{array}{l}
\bar x = \tilde M + \tilde c x -  y^2, \\
\displaystyle \bar y = -\frac{\tilde M}{\tilde c} +
\frac{1}{\tilde c} y + \frac{1}{\tilde c}(\tilde M + \tilde c x -
y^2)^2,
\end{array}
\right.
 \label{frmex}
\end{equation}
in which the coordinates $(x,y)$ and the parameters $(\tilde M,
\tilde c)$ can take arbitrary values except $\tilde{c}=0$. The region
$\tilde c<0$ will stand for the ``inner'' tangency case and $\tilde c>0$ for the ``outer'' one. Its bifurcation diagram is showed in Figure~\ref{Figdiag31}.
The map~(\ref{frmex}) is, in fact, the
product of two H\'enon maps with Jacobian $-\tilde c$ and $-\tilde c^{-1}$,
(see equations~(\ref{frmH12})).
Thus, we can state (see also~\cite{T03b}) that map (\ref{frmex})
has a complicated dynamics in the corresponding parameter intervals.
The latter  means, in particular, that all fixed points become saddles
and all of them have homoclinic and heteroclinic intersections for all values of
the parameter $\mu$ including (quadratic) tangencies for dense subsets --
\emph{Newhouse phenomenon}.

An analogous ``homoclinic tangle'' can be observed for  the
first-return map $T_{km}$ (see Lemma~\ref{lm:frm}). However,
although the map~(\ref{frmex}) is reversible and conservative (its
Jacobian is identically $1$), the original first-return map
$T_{km}$ is also reversible but not conservative in general (see
Lemma~\ref{lmnsym}). Precisely, it will be shown that in some
regions of the space of parameters $(\tilde{c},\tilde{M})$ the map
$T_{km}$ possesses chaotic dynamics and has four saddle fixed
points, two of them symmetric conservative and a symmetric couple
of fixed points (that is, symmetric one to each other and with
Jacobian greater and less than $1$, respectively). According
to~\cite{GST97,LSt,GST07}, the following result holds:

\begin{theorem} Let $\left\{ f_\mu \right\}$ be the one-parameter family of reversible maps from Theorem~\ref{th:main1}. Then, in any
segment  $[-\varepsilon,\varepsilon]$ of values of $\mu$, there are Newhouse intervals with mixed dynamics connected with an abundance of attracting, repelling  and elliptic periodic orbits. This is, values of parameters corresponding to maps $f_{\mu}$ exhibiting simultaneously infinitely many periodic orbits of
all these types form a residual set (of second category) in these intervals.
\label{th:col-new}
\end{theorem}

\section{Proof of Theorem~\ref{th:main1}}
\label{sec:prTh1}

\subsection{Preliminary geometric and analytic constructions}
\label{sec:pelim}

To ease the reading all the proofs of the lemmas of this section have been deferred to Sections~\ref{se:crossL1L2},~\ref{se:technical_lemmas} and~\ref{se:birkhoff}.

\begin{figure}[htb]
\begin{center}
\includegraphics[width=10cm]{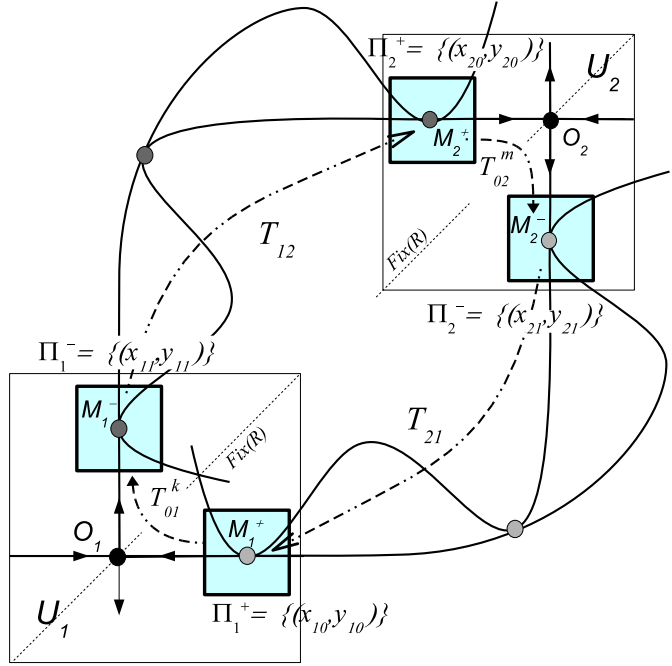} 
\caption{Schematic actions of the local ($T_{01}^k$
and $T_{02}^m$) and global ($T_{12}$ and $T_{21}$) maps in the
neighbourhood $U$ of the contour
$C=\{O_1,O_2,\Gamma_{12},\Gamma_{21}\}$. }
\label{fig:con2}
\end{center}
\end{figure}

Let us consider first the map $f_0$ and let $M_1^-\in U_1$, $M_2^+\in U_2$ be a pair of points of the orbit $\Gamma_{12}$ and $M_2^-\in U_2$, $M_1^+\in U_1$ be a pair of points of $\Gamma_{21}$. Consider $\Pi_i^+\subset U_i$  and $\Pi_i^-\subset U_i$ small neighbourhoods of the heteroclinic points $M_i^+$ and $M_i^-$ (see Figure~\ref{fig:con2}). Let us  assume that $(i)$ the heteroclinic points
are symmetric under the involution $R$ , i.e. $M_1^-= R(M_1^+)$ and $M_2^-=R(M_2^+)$, and $(ii)$ they are the ``last'' points on $U_1$ and $U_2$, that is, $f_0(M_i^-)\notin U_i$ (and, thus, $f_0^{-1}(M_i^+)\notin U_i$). \footnote{One can always take all neighbourhoods
$U_i,\Pi_i^-,\Pi_i^+,i=1,2$ to be also $R$-symmetric (that is $R(U_i)=U_i, R(\Pi_i^{\pm}) = \Pi_i^{\pm})$.}
Let $q$ be such an positive integer that $M_2^+=f_0^q(M_1^-)$ (and, thus, $M_1^+=f_0^q(M_2^-)$).

Consider now the map $f_\mu$. Denote $T_{0i} \equiv {f_\mu}\bigl|_{U_i}$, $i=1,2$. The maps
$T_{01}$ and $T_{02}$ are called {\em the local maps}. We introduce also the so-called {\it global
maps} $T_{12}$ and $T_{21}$ by the following relations: $T_{12}\equiv f^q_\mu :\Pi_1^-\to\Pi_2^+$
and $T_{12}\equiv f^q_\mu :\Pi_2^-\to\Pi_1^+$ (see Figure~\ref{fig:con2}). Then the {\em
first-return map} $T_{km}:\Pi_1^+\mapsto\Pi_1^+$ is defined by the following composition of maps and neighbourhoods:
\begin{equation}
\Pi_1^+ \;\stackrel{T_{01}^k}{\longrightarrow}\;\;\Pi_1^- \;\stackrel{T_{12}}{\longrightarrow}\;\; \Pi_2^+ \;\stackrel{T_{02}^m}{\longrightarrow}\;\;\Pi_2^-\;\stackrel{T_{21}}{\longrightarrow}\;\; \Pi_1^+
\label{TkmPi}
\end{equation}
Denote local coordinates on $\Pi_i^+$ and $\Pi_i^-$ as  $(x_{0i},y_{0i})$ and $(x_{1i},y_{1i})$,
respectively. Then the chain (\ref{TkmPi}) can be represented (in coordinates) as
\[
(x_{01},y_{01})
\stackrel{T_{01}^k}{\longmapsto}(x_{11},y_{11})
\stackrel{T_{12}}{\longmapsto}(x_{02},y_{02})
\stackrel{T_{02}^m}{\longmapsto}(x_{12},y_{12})\stackrel{T_{21}}{\longmapsto}
(\bar x_{01},\bar y_{01}).
\]

As usually, we need such local coordinates on $U_1$ and $U_2$ in which the maps $T_{01}$ and
$T_{02}$ have their simplest form.
We can not assume the maps $T_{0i}$ are
linear, since by condition {\sf [A]}, only
$C^1$-linearisation is possible here.
Therefore, we consider
such $C^{r-1}$-coordinates in which the local maps have the so-called {\em main normal
form} or normal form of the first order.

\begin{lm}[Main normal form of a saddle map]
Let a $C^r$-smooth map $T_0$ be reversible with
$\dim\mbox{Fix}(T_0) = 1$. Suppose that $T_0$ has a saddle fixed
(periodic) point $O$ belonging to the line
$\mbox{Fix}(T_0)$ and
having multipliers
$\lambda$ and $\lambda^{-1}$, with $|\lambda|<1$. Then
there exist $C^{r-1}$-smooth local coordinates near $O$ in
which the map $T_0$ (or $T_0^n$, where $n$ is the period of $O$)
can be written in the following form:
\begin{equation}
T_0:\;\;
\begin{array}{l}
\bar x = \lambda x (1 + h_1(x,y)xy)\\
\bar y = \lambda^{-1} y(1 + h_2(x,y)xy),
\label{GoodFormSaddle}
\end{array}
\end{equation}
where $h_1(0)=- h_2(0)$. The map (\ref{GoodFormSaddle}) is reversible with respect to the standard linear involution $(x,y) \mapsto (y,x)$.
\label{LmSaddle}
\end{lm}

When proving this lemma one can deduce more ``descriptive'' properties of this map. Namely, one can show that it could be written in the so-called {\em cross-form} (see Section~\ref{subse:cross_form}) as follows:
\begin{equation}
T_0:\;\;
\begin{array}{l}
\bar x = \lambda x  + \hat h(x,\bar y)x^2\bar y,\\
y = \lambda \bar y + \hat h(\bar y,x)x\bar y^2.
\label{CrossFormSaddle}
\end{array}
\end{equation}

Notice that when $T_0$ is linear, i.e. $T_0:\; \bar x = \lambda x, \bar y = \lambda^{-1}y$, one can
easyly find formulas for its iterates $T_0^j$, $j\in \Z$. Namely, one can write it either as $x_j \lambda^j x_0,  y_j = \lambda^{-j} y_0$ or, in cross-form, as $x_j = \lambda^j x_0,  y_0 = \lambda^{j} y_j$.
If $T_0$ is non-linear, then its cross-form equations exist too. In particular, the
following result holds:

\begin{lm}[Iterations of the local map]
Let $T_0$ be a saddle map written in the main normal form (\ref{GoodFormSaddle}) (or
(\ref{CrossFormSaddle})) in a small neighbourhood $V$ of $O$. Let us consider points $(x_0,y_0),\dots,(x_j,y_j)$ from $V$ such that $(x_{l+1},y_{l+1})= T_0(x_{l},y_{l})$, $l=0,\dots,j-1$. Then  one has
\begin{eqnarray}
x_j &=& \lambda^j x_0 \left(1 + j \lambda^j h_j(x_0,y_j)\right), \\
y_0 &=& \lambda^j y_j \left(1 + j
\lambda^j h_j(y_j, x_0)\right), \nonumber
\label{LocalMapk}
\end{eqnarray}
where
the functions  $h_j(y_j, x_0)$ are uniformly bounded with respect to $j$ as well all their derivatives up to order $r-2$. \label{LmLocalMap}
\end{lm}

\begin{remark}

\begin{itemize}
\item[(a)]
Both lemmas~\ref{LmSaddle} and~\ref{LmLocalMap}
remain true if $T_0$ depends on parameters. Moreover, if the
initial $T_0$ is $C^r$ with respect to coordinates and parameters,
then the normal form (\ref{GoodFormSaddle}) is $C^{r-1}$ with
respect to coordinates and $C^{r-2}$ with respect to parameters
(see~\cite{GST07}).

\item[(b)]
It follows from Lemma~\ref{LmSaddle} that the involution $L(x,y)=(y,x)$
is very convenient for the construction of symmetric saddle maps. Moreover, this involution is (locally)
smoothly equivalent to any other involution $R$ with $\dim\;\mbox{Fix}(R)=1$ (see~\cite{MZ55}). Thus, our assumption on a concrete form of the involution for the maps
$f_\mu$ does not lead to loss of generality.

\item[(c)]
Similar results related to finite-smooth normal forms of
saddle maps were established in \cite{GS90,GS00,GG04,GST07} for
general, near-conservative and conservative maps. In this paper
we, in fact, modify the corresponding proofs adapting them to the
reversible case.
\end{itemize}
\end{remark}

\subsection{Construction of the local and global maps}
\label{sec:locglob}

\begin{figure}[ht]
\begin{tabular}{cc}
\psfig{file=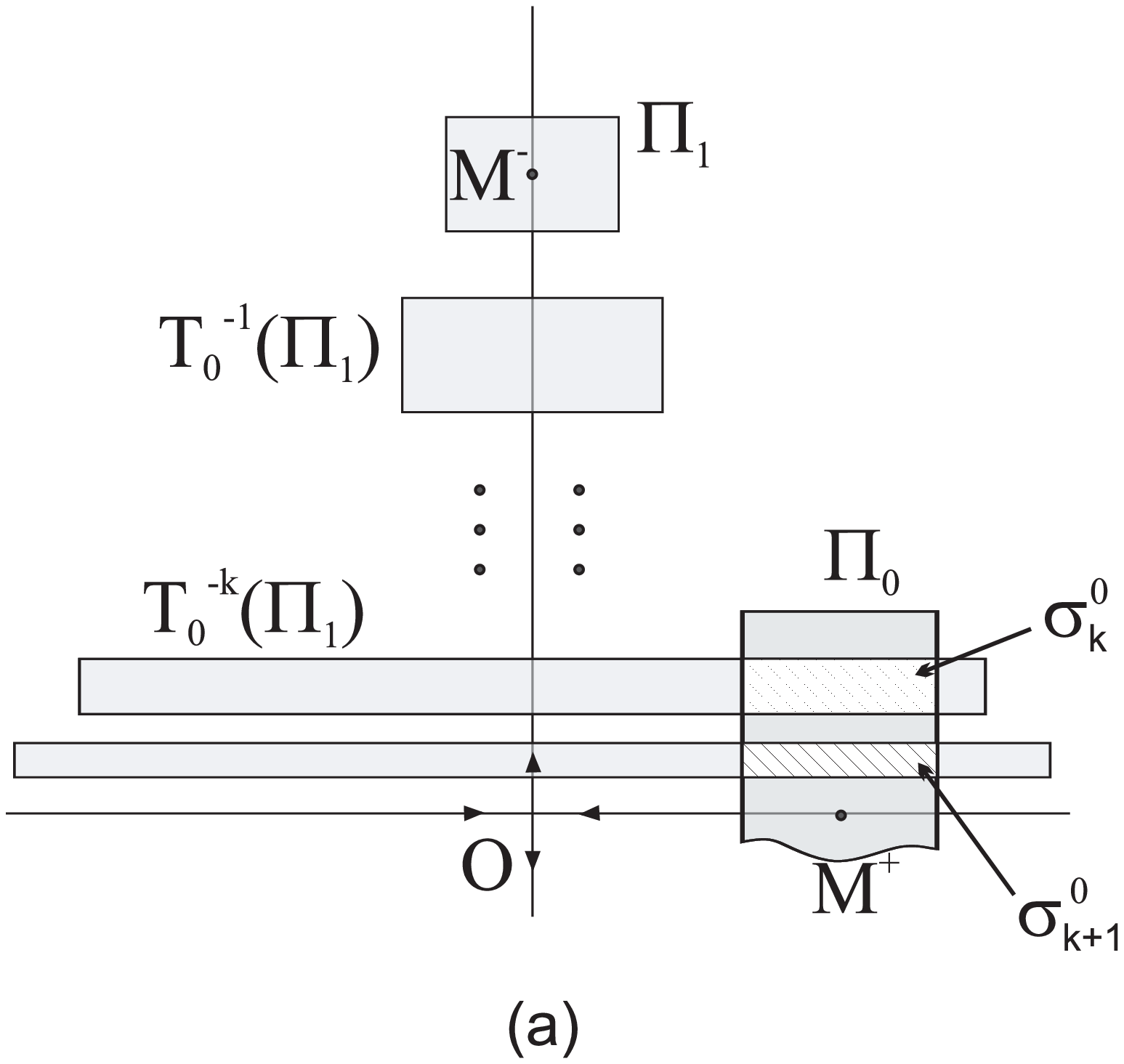,width=55mm,angle=0} \hspace{0.2cm}
\psfig{file=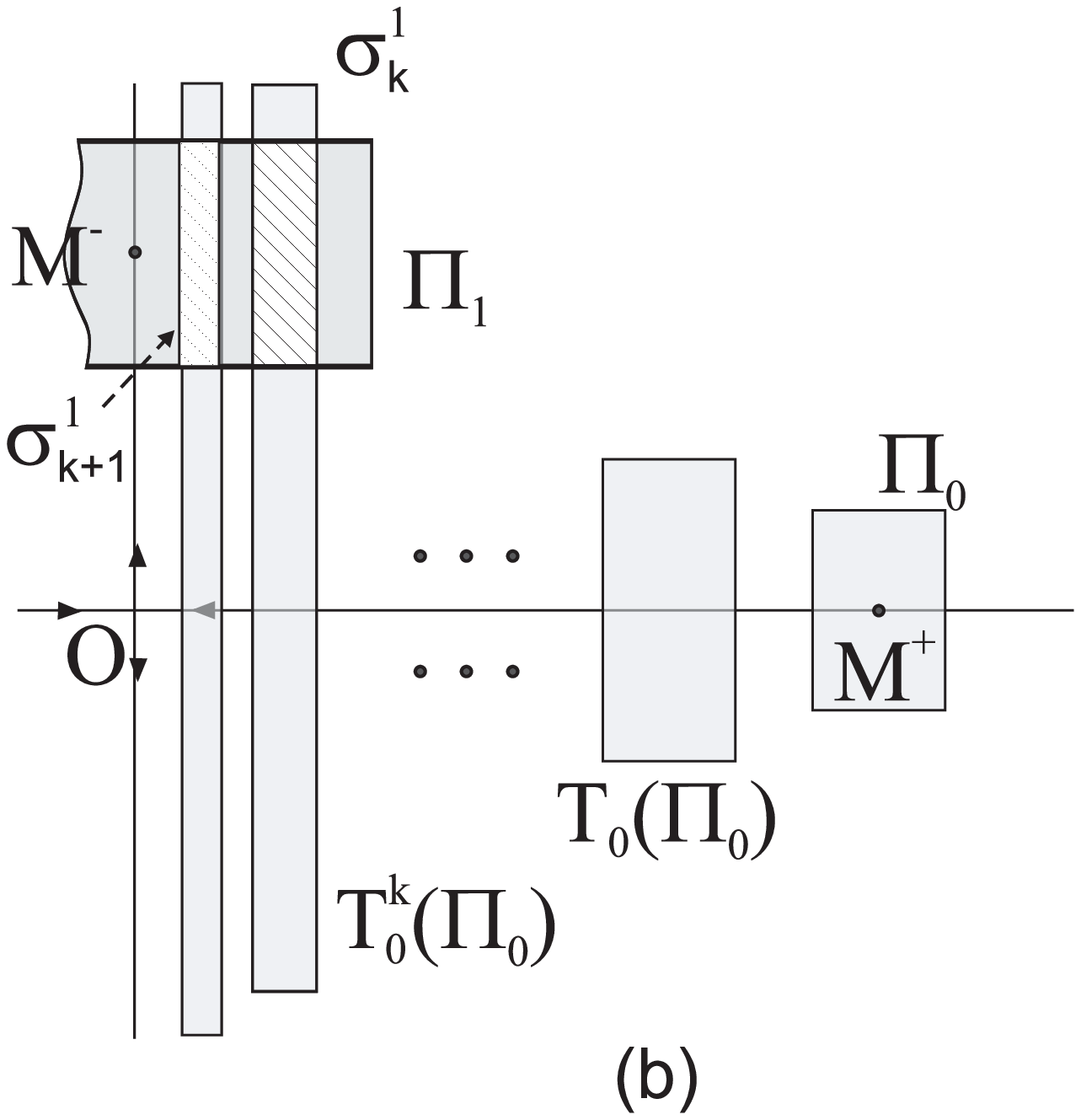,width=50mm,angle=0}
\end{tabular}
\caption{
{\footnotesize A geometry of creation of both domains of definition $\sigma_i^0\subset\Pi^+$  (a)
and domains of the range $\sigma_i^1\subset\Pi^-$ (b) for the maps $T_0^i:\Pi^+\to\Pi^-$. }}
\label{fig:strips}
\end{figure}

By Lemma~\ref{LmLocalMap}, we can choose in $U_1$ and $U_2$ local coordinates  $(x_1,y_1)$ and
$(x_2,y_2)$, respectively, such that the maps $T_{01}$ and $T_{02}$ take the following form:
\[
T_{01}:\;\; \bar x_1 = \lambda_1 x_1  + h_1^1(x_1,y_1)x_1^2y_1,\quad \bar y_1 = \lambda_1^{-1} y_1 +
h_2^1(x_1,y_1)x_1y_1^2,
\]
and
\[
T_{02}:\;\; \bar x_2 = \lambda_2 x_2  + h_1^2(x_2,y_2)x_2^2y_2,\quad \bar y_2 = \lambda_2^{-1} y_2 +
h_2^2(x_2,y_2)x_2y_2^2.
\]
Furthermore, in these coordinates, the local stable and unstable invariant manifolds of both points $O_1$ and $O_2$ are straightened: $x_i=0$ is the equation of $W^u_{loc}(O_i)$  and $y_i=0$ is the equation of
$W^s_{loc}(O_i)$, $i=1,2$. Then, we can write the $(x,y)$-coordinates of the chosen
heteroclinic points as follows: $M_1^+(x_1^+,0)$, $M_1^-(0,y_1^-)$, $M_2^+(x_2^+,0)$ and
$M_2^-(0,y_2^-)$. Besides, because of the reversibility, we have that
\begin{equation}
x_1^+=y_1^- = \alpha_1^*,\qquad x_2^+=y_2^- = \alpha_2^*
\label{xeqy}
\end{equation}
We assume that $T_{0i}(\Pi_i^{+})\cap\Pi_i^{+}=\emptyset$ and
$T_{0i}^{-1}(\Pi_i^{-})\cap\Pi_i^{-}=\emptyset$, $i=1,2$.
Then the domain of definition of the successor map from $\Pi_i^{+}$
into $\Pi_i^-$ under iterations of $T_{0i}$ consists
of infinitely many
non-intersecting strips $\sigma_{j}^{0i}$ which belong to $\Pi_i^+$ and accumulate at
$W^s_{loc}(O_i)\cap\Pi_i^+$ as $j\to\infty$. Analogously, the range of this map consists of infinitely many  strips $\sigma_{j}^{1i}= T_{0i}^j(\sigma_j^{0i})$ belonging to $\Pi_i^-$ and accumulating at $W^u_{loc}(O_i)\cap\Pi_i^-$ as $j\to\infty$ (see Figure~\ref{fig:strips}).

It follows from Lemma~\ref{LmLocalMap} that the map $T_{01}^k:\sigma_k^{01}\mapsto\sigma_k^{11}$
can be written in the following form (for large enough values of $k$)
\begin{equation}
T_{01}^k\;:\; \left\{
\begin{array}{l}
x_{11} = \lambda_1^k x_{01} (1 + k\lambda_1^k h_k^1(x_{01},y_{11}) )\\
y_{01} = \lambda_1^k y_{11} (1 + k\lambda_1^k h_k^1(y_{11},x_{01}) )
\end{array}
\right. \label{T1k}
\end{equation}
and an analogous formula takes place for the map $T_{02}^m:\sigma_m^{02}\mapsto\sigma_m^{12}$:
\begin{equation*}
T_{02}^m :\; \left\{\begin{array}{l}
x_{12}= \lambda_2^m x_{02}(1 + m\lambda_2^m h_m^2(x_{02},y_{12}) )\\
y_{02} = \lambda_2^m y_{12}(1 + m\lambda_2^m h_m^2(y_{12},x_{02}) )
\end{array}
\right. \label{T2m}
\end{equation*}

We write now the global map $T_{12}$  in the following form
\begin{equation}
\fl T_{12}
\left\{
\begin{array}{l}
x_{02} - x_2^+ = F_{12}(x_{11},y_{11}-y_1^-,\mu) \equiv a x_{11} + b (y_{11}-y_1^-) + \\
\qquad\qquad l_{02}(y_{11}-y_1^-)^2 + \varphi_1(x_{11},y_{11},\mu),\\
y_{02} = G_{12}(x_{11},y_{11}-y_1^-,\mu) \equiv  \mu + c x_{11} +
d (y_{11}-y_1^-)^2  + \\
\qquad\qquad
f_{11}x_{11}(y_{11}-y_1^-)+f_{03}(y_{11}-y_1^-)^3+ \\
\qquad\qquad \varphi_2(x_{11},y_{11},\mu),
\end{array}
\right.
\label{T12}
\end{equation}
where $F_{12}(0)=G_{12}(0)=0$  since $T_{12}(M_1^-) = M_2^+$ at
$\mu=0$ and
\begin{equation*}
\begin{array}{l}
\varphi_1=O(|y_{11}-y_1^-|^3) + x_{11} O(\|(x_{11},y_{11}-y_1^-)\|),\\
\varphi_2=O(|x_{11}|^2)+ O(|y_{11}-y_1^-|^4) +
O(x_{11}(y_{11}-y_1^-)^2).
\end{array}
\label{T12hot}
\end{equation*}
Since the curves
$T_{12}\left(W^u_{loc}(O_1):\{x_{11}=0\}\right)$ and
$W^s_{loc}(O_2):\{y_{02}=0\}$ have a quadratic tangency at
$\mu=0$, it implies that
\[
\frac{\partial G_{12}(0)}{\partial y_{11}} = 0,\;\;
\frac{\partial^2 G_{12}(0)}{\partial y_{11}^2} = 2d \neq 0.
\]
The Jacobian $J(T_{12})$ has, obviously, the following form:
\[
J(T_{12})= -bc + af_{11} x_{11} + Q\; (y_{11}-y_1^-)+ O\left(x_{11}^2+
(y_{11}-y_1^-)^2\right),
\]
where
\begin{equation}
\label{RFRMJac}
Q = 2ad- bf_{11} - 2c l_{02}.
\end{equation}
Now condition \textsf{[C]} can be formulated more precisely. Namely, we
require that\footnote{Notice that the analogous phenomenon (of
influence of high order terms on dynamics) was discovered in
\cite{GG00,GG04} when studying bifurcations homoclinic tangencies
to a saddle fixed point with the unit Jacobian.}
\begin{equation*}
Q \;=\; \frac{\partial J(T_{12})}{\partial
y_{11}}{\Bigl|_{{(x_{11}=0,y_{11}=y_1^-,\mu=0)}}}\;\neq\; 0
\label{Qi0}
\end{equation*}

Concerning the global map $T_{21}$, we cannot write it now in an
arbitrary form. The point is that after written a
formula for the map $T_{12}$ it is necessary to use the reversibility relations
to get the one associated to it:
\begin{equation*}
T_{21} = R\; T_{12}^{-1}\;R^{-1}, \qquad  T_{12} = R\; T_{21}^{-1}\;R^{-1}
\label{T12-T21}
\end{equation*}
for constructing $T_{21}$. Then, by (\ref{T12}), we obtain that the map
$T_{21}^{-1}:\Pi_1^+\{(x_{01},y_{01})\}\mapsto \Pi_2^-\{(x_{12},y_{12})\}$ must be written as follows
\begin{equation}
\fl T_{21}^{-1} \left\{
\begin{array}{l}
x_{12} = G_{12}(y_{01},x_{01}-y_1^-,\mu) = \\
\quad \mu + c y_{01} + d
(x_{01}-y_1^-)^2 + f_{11}
y_{01}(x_{01}-y_1^-) + \\
\quad f_{13}(x_{01}-y_1^-)^2 + \varphi_2(y_{01},x_{01},\mu),\\
y_{12} - x_2^+ = F_{12}(y_{01},x_{01}-y_1^-,\mu) =  \\
\quad  a y_{01} +
b (x_{01}-y_1^-) + l_{02}(x_{01}-y_1^-)^2 + \varphi_1(y_{01},x_{01},\mu)\\
\end{array}
\right.
 \label{T21-}
\end{equation}
Relation (\ref{T21-}) allows to define the map
$T_{21}:\Pi_2^-\{(x_{12},y_{12})\}\mapsto\Pi_1^+\{(x_{01},y_{01})\}$, but in implicit
form.

\subsection{Construction of the first-return maps $T_{km}$ and the Rescaling Lemma}
\label{sec:resclm}

Now, using relations (\ref{T1k})--(\ref{T21-}), we can construct
the first-return map $T_{km}= T_{21}T_{02}^mT_{12}T_{01}^k$
defined on
the strip
$\sigma_k^{01}\subset\Pi^+_1$.
Recall that any fixed point of $T_{km}$ corresponds to a
single-round periodic orbit of type $(k,m)$ of period $(k+m+2q)$.
However, we do not state the problem of studying the maps $T_{km}$
for all large $k$ and $m$.  We suppose $k$ and $m$ are large
enough integers such that
\begin{equation}
\lambda_1^k \simeq \lambda_2^m.
\label{keqm}
\end{equation}
In other words, both values of $\lambda_1^k\lambda_2^{-m}$ and $\lambda_1^{-k}\lambda_2^{m}$ are
{\em uniformly} separated from $0$ and $\infty$ as $k,m\to\infty$. Then the following result holds.

\begin{lm}[The rescaling lemma]
%
Let the map $f_0$  satisfy  conditions {\sf [A]-[B]} and $f_{\mu}$
be a general unfolding in the class of reversible maps. Suppose
$k$ and $m$ are large enough integer numbers satisfying
relation~(\ref{keqm}). Then one can introduce coordinates (called
``rescaled coordinates'') in such a way that the first-return map
$T_{km}$ takes the form
\begin{equation}
\begin{array}{l}
M + c \bar y + d \bar x^2 + f_{11}\lambda_1^{k} \bar x\bar y+ f_{03} \lambda_1^{k}\bar x^3 = \\
\qquad = b \lambda_2^m \lambda_1^{-k} y + a \lambda_2^m x +
l_{02}\lambda_2^{m}y^2 + O(k\lambda_1^{2k}),\\
\\
M + c x + d y^2 + f_{11}\lambda_1^{k} xy+ f_{03}\lambda_1^{k}y^3 = \\
\qquad b \lambda_2^m \lambda_1^{-k} \bar x + a \lambda_2^m \bar y + l_{02}\lambda_2^{m}\bar x^2 +
O(k\lambda_1^{2k}),
\end{array}
\label{frmSO}
\end{equation}
where
\begin{equation}
M=\lambda_1^{-2k}\left(\mu + c \lambda_1^k \alpha_1^*(1+\dots)  - \lambda_2^m \alpha_2^*(1 +
\dots)\right)
\label{Mmu}
\end{equation}
and ``$\ldots$'' stands for some coefficients tending to zero as
$k,m\to\infty$. Notice that the domain of definition of the new
coordinates $x\sim\lambda_1^{-k}(x_{01}-\alpha_1^*)$,
$y\sim\lambda_1^{-k}(y_{11}-\alpha_1^*)$ and parameter $M$ cover
all finite values as $k,m \to \infty$.
\label{lm:frm}
\end{lm}

\subsection{On bifurcations of fixed points of the first-return maps $T_{km}$}
\label{sec4}

We study bifurcations in the first-return map $T_{km}$ using its rescaled form (\ref{frmSO}).
If we neglect
in~(\ref{frmSO}) all asymptotically small terms (as $k,m\to\infty$),
we obtain the following truncated form for $T_{km}$
\begin{equation}
\begin{array}{l}
M + c \bar y + d \bar x^2 =  \beta_{km} y, \;\;
\beta_{km}\bar x = M + c x + d y^2,
\end{array}
\label{frmFO}
\end{equation}
where $\beta_{km}=b \lambda_1^{-k }\lambda_2^m$. Rescale the
coordinates
$$
x\;=\; -\frac{\beta_{km}}{d} x_{new},\;\; y\;=\;
-\frac{\beta_{km}}{d} y_{new}.
$$
Then map (\ref{frmFO}) is rewritten in the following form, that we denote with $H$:
\begin{equation}
H: \ \
\tilde M + \tilde c \bar y - \bar x^2 =   y, \qquad
\bar x = \tilde M + \tilde c x -  y^2,
\label{stfrmFO}
\end{equation}
where
\begin{equation}
\tilde M = -\frac{d}{\beta_{km}^2}\;M,\;\;  \tilde c = \frac{1}{\beta_{km}}\equiv
\frac{c}{b}\lambda_1^k\lambda_2^{-m}\;
\label{frmMC}
\end{equation}
Notice that $H$ depends on two
parameters $\tilde M$ and $\tilde c$ which can take arbitrary
values, except for $\tilde c =0$ (according to conditions $bc\neq
0$ and $0<|\lambda_i|<0$). Thus, two mainly different scenarios take place: with $\tilde c<0$ and $\tilde c>0$.

Observe that the map $H$ can be expressed in the explicit form (\ref{frmex}). Moreover, it can be
represented as the superposition $H=\mathcal{H}_2\circ \mathcal{H}_1$ of two quadratic (H\'enon) maps
\begin{equation*}
\begin{array}{l}
\mathcal{H}_1:\{\;\; x_1 = y_0, \qquad y_1 =  \tilde M + \tilde c x_0 -  y_0^2, \\
\displaystyle \mathcal{H}_2:\{\;\; x_2 = y_1,\qquad y_2 = -\frac{\tilde M}{\tilde c} +
\frac{1}{\tilde c} x_1 + \frac{1}{\tilde c} y_1^2.
\end{array}
\label{frmH12}
\end{equation*}
The Jacobians of these maps are constant and inverse: $J(\mathcal{H}_1)= -\tilde c$ and $J(\mathcal{H}_2)=-\tilde
c^{-1}$. Therefore, the resulting map $H=\mathcal{H}_2\circ \mathcal{H}_1$ is a quadratic map with  the Jacobian equal to
$1$ and so area-preserving.\footnote{Notice that a non-trivial technique proposed in
\cite{T03b} based on considering superpositions of H\'enon-like maps, allows to
deduce a series of quite delicious generic properties demonstrating richness of chaos in
non-hyperbolic area-preserving maps.}

\begin{figure}[ht]
\begin{center}
\includegraphics[width=10cm]{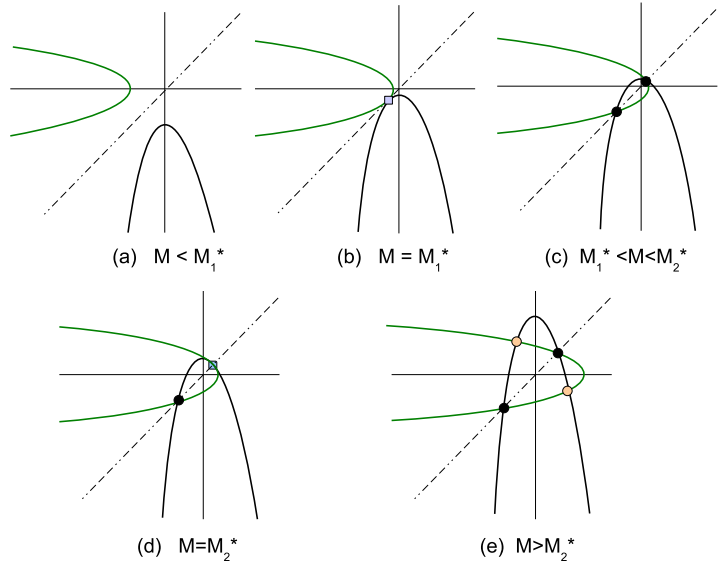} 
\caption{A geometric interpretation for a structure
of fixed points of map (\ref{stfrmFO}) as intersection points of
two symmetric parabolas.}
 \label{Fig4}
\end{center}
\end{figure}

Form (\ref{stfrmFO}) of $H$ allows to give a rather simple geometric interpretation of the
bifurcations of fixed points. The coordinates $(x,y)$ of these fixed points must satisfy
the equations
\begin{equation}
y(1-\tilde c)= \tilde M - x^2, \qquad  x(1-\tilde c)= \tilde
M - y^2
\label{geomH12}
\end{equation}
Let us hold fixed $\tilde c$ and suppose that $\tilde c \neq 1$. Then equations of (\ref{geomH12})
define on the $(x,y)$-plane two parabolas which are symmetric with respect to the bisectrix $y=x$.
Intersection points of the parabolas are also fixed points of $H$. When $\tilde M$ varies the
parabolas ``move'' and, as a result, the number of intersection points can change (i.e. bifurcations in
$H$ occur). See Figure~\ref{Fig4} in which the case $\tilde c <1$ is illustrated: (a) the parabolas
do not intersect if $\tilde M<M_1^*\equiv -\frac{1}{4}(\tilde c-1)^2$; (b) the parabolas are
touched (quadratically) to the bisectrix and one to other,  if $\tilde M=M_1^*$; (c) they have two
(symmetric) intersection points if $M_1^*<\tilde M<M_2^*\equiv \frac{3}{4}(\tilde c-1)^2$; (d) they
have a cubic (symmetric) tangency when $\tilde M= M_2^*$ and, finally, (e) the parabolas have four
intersection points (two symmetric points and a symmetric couple of points) if $\tilde M>M_2^*$.
An analogous picture takes place for the case $\tilde c >1$ (the parabolas have their branches in the opposite directions). The
case $\tilde c =1$ is very special. Here, the equation
(\ref{geomH12}) takes the form $0= \tilde M - x^2, 0 = \tilde M -
y^2$ and, thus, a certain $"0-4"$-bifurcation occurs at $\tilde
M=0$: the map $H$ has no fixed points for $\tilde M<0$ and 4 fixed
points appear immediately when $\tilde M$ becomes positive.

\begin{figure} [bht]
\begin{center}
\includegraphics[width=12cm]{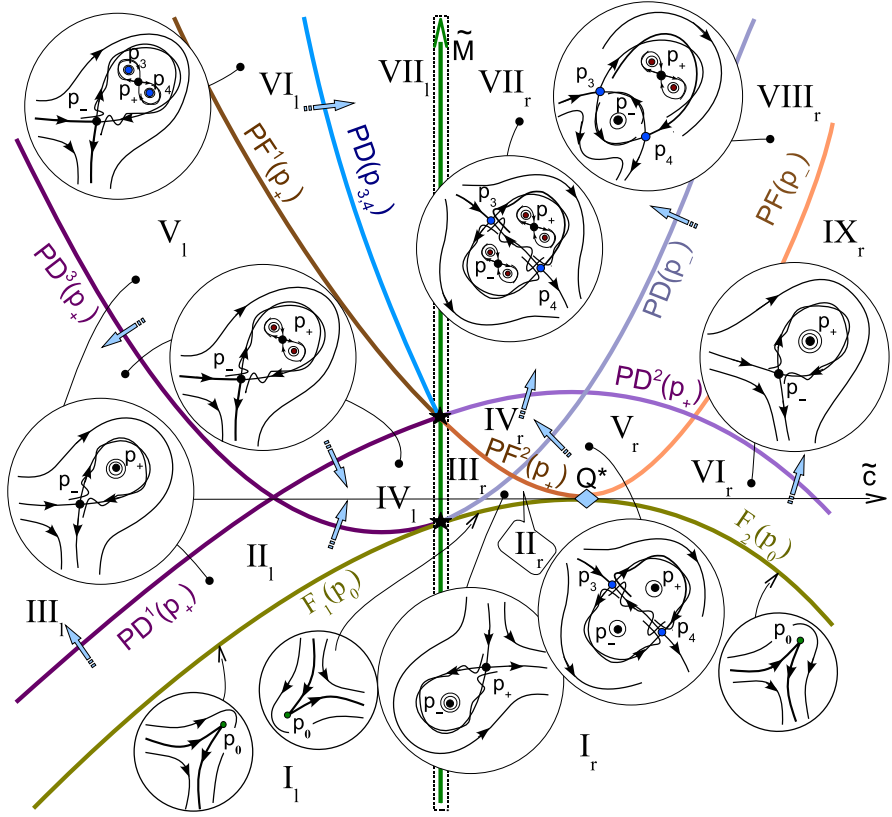} 
\caption{Elements of the bifurcation diagram for
the map $H$.}
 \label{Figdiag31}
\end{center}
\end{figure}

More details concerning bifurcations of fixed points of $H$ are illustrated in
Figure~\ref{Figdiag31} where
principal elements of the bifurcation diagram on the $(\tilde
c,\tilde M)$-plane are represented.
Notice that the line $\tilde c= 0$ is singular and, therefore,
there are no transitions between the half-planes $\tilde c <0$ and $\tilde
c >0$. In particular, this means that bifurcation
curves must ``terminate'' on the line $\tilde c = 0$ (two such
terminated points are denoted in Figure~\ref{Figdiag31} as black
stars). Besides, three types of bifurcation curves are represented
in the figure: fold ($F$), period-doubling ($PD$) and pitch-fork
($PF$).

The curves $\mbox{F}_1$ and $\mbox{F}_2$ having the same equation
\begin{equation}
\tilde M=-\frac{1}{4}(\tilde c-1)^2
\label{bfold12}
\end{equation}
but with $\tilde c <0$ and $\tilde c >0$, respectively, relate to
a  conservative fold-bifurcation. If $\tilde M <
-\frac{1}{4}(\tilde c-1)^2$, i.e. $(\tilde c,\tilde M)\in I_l\cup
I_r$,  the map $H$ has no fixed points; if $\tilde M >
-\frac{1}{4}(\tilde c-1)^2$, the map $H$ has two symmetric fixed
points $P^+=(p_+,p_+)$ and $P^-=(p_{-},p_{-})$, where
\begin{equation}
p^{\pm}=\frac{\tilde c-1 \pm \sqrt{(\tilde c - 1)^2+4\tilde
M}}{2}. \tilde M=-\frac{1}{4}(\tilde c-1)^2
\label{p+p-}
\end{equation}
Notice that the point $Q^* = (\tilde c=1,\tilde M =0)\in
\mbox{F}_2$ corresponds to a degenerate fold-bifurcation:
simultaneously four fixed points, two symmetric and a symmetric
couple, are born at the transition $\mbox{I}_r\rightarrow
\mbox{V}_r$.

In the case $\tilde c<0$, the point $p_-$ is always a saddle and
it does not bifurcate any more. This is not the case of the point $p_+$
which can undergo both period-doubling and pitch-fork
bifurcations. The period-doubling bifurcation curves
\begin{equation*}
\begin{array}{l}
\displaystyle PD^1(p_+):\;\; \tilde M= 1 -\frac{1}{4}(\tilde
c-1)^2, \;\; \tilde
c <0; \\
\displaystyle PD^2(p_+):\;\; \tilde M= \frac{(c+1)(3c-1)}{4}, \;\;
\tilde c <0.
\end{array}
\label{bpd+}
\end{equation*}
are represented in Figure~\ref{Figdiag31} by ``grey arrows'' which
indicate directions of birth of period-$2$ points. The curve
\begin{equation*}
\begin{array}{l}
PF^1(p_+):\;\; \tilde M= \frac{3}{4}(\tilde c-1)^2, \;\; \tilde
c <0 \\
\end{array}
\label{bpfor2}
\end{equation*}
relates to the pitch-fork bifurcation: when crossing this curve
(in the direction $\mbox{V}_l\rightarrow \mbox{VI}_l$) the point
$p_+$ becomes saddle and two asymmetric {\em elliptic} fixed
points $p_3$ and $p_4$ are born in its neighbourhood. The point
$p_+$ does not bifurcate any more, whereas, the points $p_3$ and $p_4$
undergo simultaneously period-doubling bifurcation at crossing the
curve
\begin{equation*}
PD(p_{3,4}):\;\; \tilde M=\frac{(1-3\tilde c)(3-\tilde c)}{4}, \quad
\tilde c <0
\label{bpdoubl}
\end{equation*}
Further variation of parameters in the domain $\mbox{VII}_l$, will lead to a cascade of (conservative) period-doubling bifurcations of asymmetric periodic points.

As it can be seen in Figure~\ref{Figdiag31}, the character of the
bifurcations in the case
$\tilde c>0$ is different to the one for $\tilde c <0$.
In this case, $\tilde{c}>0$, both symmetric fixed points $p_+$ and $p_-$
undergo pitch-fork and period-doubling bifurcations. The
corresponding bifurcation curves are
\begin{equation*}
\begin{array}{l}
\displaystyle PF^2(p_+): \quad \tilde M= \frac{3}{4}(\tilde c-1)^2,
\;\; 0<
\tilde c <1; \\
\displaystyle PD^3(p_+):\quad \tilde M= 1- \frac{1}{4}(\tilde
c-1)^2, \quad \tilde c >0; \\
\end{array}
\label{bpd+1}
\end{equation*}
 for the point $p_+$ and
\begin{equation}
\begin{array}{l}
\displaystyle PF(p_-): \quad \tilde M= \frac{3}{4}(\tilde c-1)^2,
\quad \tilde c >1; \\[1.2ex]
\displaystyle PD(p_-): \quad  \tilde M= \frac{(c+1)(3c-1)}{4}, \quad \tilde c >0 \\
\end{array}
\label{bpd+2}
\end{equation}
for the point $p_-$. Notice that, when crossing the curve
$PF^2(p_+)\cup PF(p_-)$, a symmetric couple of \emph{saddle} fixed
points $p_3$ and $p_4$ are born and they do not bifurcate any more
(for values of parameters $\tilde c>0$ and $\tilde M$ in the
domain lying above the curve $PF^2(p_+)\cup PF(p_-)$). One can
expect, however, that bifurcations of the symmetric fixed points
$p_+$ and $p_-$ give rise to cascades of period-doubling
bifurcations.

It should be noted that, despite the reversibility, the pitch-fork
and period-doubling bifurcations can have in $T_{km}$ a different
character in comparison with the one of the truncated map $H$. Of
course, the pitch-fork bifurcation in $T_{km}$ leads again to the
appearance of two non-symmetric fixed points $p_3$ and $p_4$ but
these points can be non-conservative. In fact, this is a general
property (it holds for open and dense set of systems) that the
following lemma shows.

\begin{lm}[Non-conservative fixed points]
The non-symmetric fixed points $p_3$ and $p_4$ of map {\rm
(\ref{frmSO})}, with $k$ and $m$ satisfying
{\rm (\ref{keqm})}, have Jacobian $J_{ns}$ and $J_{ns}^{-1}$, respectively,
with
\begin{equation}
J_{ns}=1 + \frac{Q(\eta^*-\xi^*)}{bc}\lambda_1^k + o(\lambda_1^k),
\label{FRMJac}
\end{equation}
where $Q$ is the coefficient given by formula {\rm(\ref{RFRMJac})} and $\xi^*$ and $\eta^*$ are, respectively, the $x$- and $y$-coordinate of the fixed point.
\label{lmnsym}
\end{lm}

Due to the reversibility, the fold bifurcation in the first-return
map $T_{km}$ has the same character as in the  truncated map $H$
and leads, therefore, to the appearance of two symmetric fixed
points, $p_+$ and $p_{-}$, saddle and elliptic ones. Concerning
the symmetric elliptic fixed points, we have the following result.


\begin{lm}[Symmetric elliptic
fixed points]
The point $p_+$ (resp. the point $p_-$) is generic elliptic, that
is, it is KAM-stable, for open and dense sets of values of the
parameters $(\tilde c,\tilde M)$ in the domains
$\mbox{II}_l\cup\mbox{V}_l$ and
$\mbox{IV}_r\cup\mbox{V}_r\cup\mbox{VI}_r$ (resp. in the domain
$\mbox{II}_r\cup\mbox{V}_r\cup\mbox{VIII}_r$).
\label{lmgenel}
\end{lm}

The proofs of lemmas~\ref{lm:frm}--\ref{lmgenel} are given in
Section~\ref{se:birkhoff}.

\subsection{End of the proof of Theorem~\ref{th:main1}.}

If $k$ and $m$ are large enough and having in mind~(\ref{keqm}), (\ref{Mmu}) and~(\ref{frmMC}), the following
relation between the parameters $\mu$ and $\tilde M$ holds:
\begin{equation*}
\mu = \lambda_2^m\alpha_2^*(1+\rho_k^1)-
c\lambda_1^k\alpha_1^*(1+\rho_k^2) - \frac{b^2}{d}\tilde M
\lambda_2^{2m}(1+\rho_k^3),
\label{mu12}
\end{equation*}
where $\rho_k^i,i=1,2,3,$ are some small coefficients
($\rho_k^i\to 0$ as $k\to\infty$). Using formulas
(\ref{bfold12})--(\ref{bpd+2}) for the bifurcation curves of the
truncated map (\ref{stfrmFO}), asymptotically close to the (rescaled)
first-return map (\ref{frmSO}), we
find the following expressions for the bifurcation values of $\mu$
at the statement of Theorem~\ref{th:main1}:
a value $\mu=\mu_{\rm fold}^{(k,m)}$,  which corresponds to the fold bifurcation in
$T_{km}$,
\begin{eqnarray*}
\mu_{\rm fold}^{(k,m)} &=& \lambda_2^m\alpha_2^*(1+\rho_k^1)-
c\lambda_1^k\alpha_1^*(1+\rho_k^2) + \\
&& \frac{1}{4d}(b-c\lambda_1^k\lambda_2^m)^2
\lambda_2^{2m}(1+\rho_k^3)
\end{eqnarray*}
and a value $\mu=\mu_{pf}^{(k,m)}$,
associated to the pitch-fork bifurcation in $T_{km}$ (which is
not conservative if we have in mind Lemma~\ref{lmnsym}),
\begin{eqnarray*}
\mu_{pf}^{(k,m)} &=& \lambda_2^m\alpha_2^*(1+\rho_k^1)-
c\lambda_1^k\alpha_1^*(1+\rho_k^2) - \\
&& \frac{3}{4d}(b-c\lambda_1^k\lambda_2^m)^2
\lambda_2^{2m}(1+\rho_k^3).
\end{eqnarray*}
These considerations imply Theorem~\ref{th:main1}.

\section{On applied reversible maps with mixed dynamics}
\label{se:examples}

In this section we present two concrete examples of reversible systems
where Theorem~\ref{th:main1} applies and exhibiting, therefore, mixed dynamics: periodically perturbed Duffing equation and the Pikovsky-Topaj  model~\cite{PT02} for coupled rotators.

\subsection{A periodic perturbation of the Duffing Equation}
\label{subse:Duffing}
Let us consider the following system
\begin{equation}
\label{example:system} \left\{
\begin{array}{rcl}
\dot{x} &=& y, \\
\dot{y} &=& -x + x^3 +  \varepsilon \left( \alpha + \beta y \sin \omega t \right),
\end{array}
\right.
\end{equation}
where $\alpha, \beta, \omega \in \R$ and $\varepsilon$ is an small perturbation parameter. The
unperturbed system, for $\varepsilon=0$, corresponds to the so-called Duffing system (also called \emph{Anti-Duffing} for several authors).
It is Hamiltonian, with
\[
H(x,y)=\frac{y^2}{2} + V(x), \qquad \quad V(x)=\frac{x^2}{2} - \frac{x^4}{4} - \frac{1}{4}
\]
and is (time)-reversible with respect the following linear involutions $R(x,y)=(x,-y)$ and
$S(x,y)=(-x,y)$. This system has three singular points: one elliptic at $(0,0)$ and two saddles at
$(\pm 1,0)$. Moreover, these two points are connected through two (symmetric) heteroclinic orbits
$\Gamma_h^{\pm}$.

The perturbed system, for $\varepsilon\neq 0$, is still $R$-reversible but, in principle, non
necessarily Hamiltonian. This is a particular case of a more general family of $R$-reversible
perturbations
\[
\left\{
\begin{array}{rcl}
\dot{x} &=& y, \\
\dot{y} &=& -x + x^3 +  \varepsilon g(x,y,t),
\end{array}
\right.
\]
satisfying that $g(x,-y,-t)=g(x,y,t)$. It is well known that  two (symmetric) hyperbolic periodic
orbits $\gamma_{\varepsilon}^{\pm}$ appear close to the saddle points $(\pm 1,0)$. Let us denote by
$W^{u,s}(\gamma_{\varepsilon}^{\pm})$ their corresponding unstable and stable invariant manifolds,
respectively. Generically these invariant manifolds will intersect each other transversally and
will remain close to the unperturbed heteroclinic connection. The first order in $\varepsilon$
associated to their splitting, will be given by the well-known Poincar\'e-Melnikov-Arnol'd function
\[
M(t_0) = \int_{-\infty}^{+\infty} L_{F} G \left( \Gamma_h(t) \right) \, dt,
\]
where $F(x,y)=(y,-x+x^3)$,  $G(x,y)=\varepsilon (0,\alpha + \beta y \sin \omega t)$, $\Gamma_h(t)$
is any of both unperturbed heteroclinic connections $\Gamma_h^{\pm}(t)$ and $L_F(G) = (DF)G$ stands
for the Lie derivative of $G$ with respect to $F$. Simple zeroes of $M(t_0)$ provide tangent
intersections between the invariant manifolds $W^{u,s}(\gamma_{\varepsilon}^{\pm})$. This systems
constitutes a good candidate to apply our results.

For the computation of $M(t_0)$ we consider here the positive heteroclinic orbit
$\Gamma_h=\Gamma_h^{+}$ but, by symmetry, everything applies exactly for $\Gamma_h^{-}$. Thus,
\[
\Gamma_h(t)= \left( x_h(t), y_h(t) \right) = \left( x_h(t), \dot{x}_h(t) \right) = \left( \tanh
\frac{t}{\sqrt{2}},  \frac{1}{\sqrt{2}} \sech^2 \frac{t}{\sqrt{2}} \right)
\]
and
\begin{eqnarray}
\lefteqn{ M(t_0) = \int_{-\infty}^{+\infty} (DF)G|_{(x_h(t),y_h(t),t+t_0)} \, dt = } \\
&& \int_{-\infty}^{+\infty} y_h(t) \left( \alpha + \frac{\beta}{\sqrt{2}} \left( \sech^2 \frac{t}{\sqrt{2}} \right) \sin \omega(t+t_0)  \right) \,dt = \nonumber \\
&& \frac{\varepsilon \alpha}{\sqrt{2}}  \int_{-\infty}^{+\infty} \sech^2\frac{t}{\sqrt{2}} \, dt +
\frac{\varepsilon \beta}{2} \int_{-\infty}^{+\infty} \left( \sech^4 \frac{t}{\sqrt{2}} \right) \sin \omega(t+t_0) \, dt = \nonumber \\
&& \frac{\varepsilon \alpha}{\sqrt{2}} I_1 + \frac{\varepsilon \beta}{2} I_2.
\label{example:Mto_formula}
\end{eqnarray}
Concerning $I_1$ it is straightforward to check that its value is $2$. Regarding $I_2$, it is more
convenient to compute the integral
\[
\int_{-\infty}^{+\infty}  \left( \sech^4 \frac{t}{\sqrt{2}} \right) \rme^{\rmi \omega(t+t_0)} \, dt
\]
using the method of residues. Indeed, from it we can derive that
\begin{eqnarray*}
\int_{-\infty}^{+\infty} \left( \sech^4 \frac{t}{\sqrt{2}} \right) \sin \omega(t+t_0) \, dt &=&
\frac{2\pi}{3} \frac{\omega^2(\omega^2 + 2)}{\sinh \frac{\omega\pi}{\sqrt{2}}} \sin \omega t_0, \\
\int_{-\infty}^{+\infty} \left( \sech^4 \frac{t}{\sqrt{2}} \right) \cos \omega(t+t_0) \, dt &=&
\frac{2\pi}{3} \frac{\omega^2(\omega^2 + 2)}{\sinh \frac{\omega\pi}{\sqrt{2}}} \cos \omega t_0
\end{eqnarray*}
and, substituting in~(\ref{example:Mto_formula}), we get
\begin{eqnarray*}
M(t_0) &=& \varepsilon \left( \alpha\sqrt{2} + \frac{\beta\pi}{3} \frac{\omega^2(\omega^2+2)}{\sinh \frac{\omega\pi}{\sqrt{2}}} \sin\omega t_0 \right) = \\
&&\frac{3\sinh \frac{\omega \pi}{\sqrt{2}}}{\omega^2 (\omega^2+2)} \varepsilon \left( \alpha {\cal
P}(\omega) + \beta \sin\omega t_0 \right),
\end{eqnarray*}
provided we define
\[
{\cal P}(\omega)= \frac{\sqrt{2}\omega^2 (\omega^2+2)}{3 \sinh\frac{\omega\pi}{\sqrt{2}}}.
\]
Therefore, for small values of $\varepsilon$ we have: $(i)$ if $|\beta/\alpha| > {\cal P}(\omega)$
then $W^{u}(\gamma_{\varepsilon}^{-})$ and $W^{s}(\gamma_{\varepsilon}^{+})$ intersect; $(ii)$ if
$\beta/\alpha|<{\cal P}(\omega)$ they do not intersect each other and $(iii)$ if
$|\beta/\alpha|={\cal P}(\omega)$ then $M(t_0)$ has zeroes which are double but not triple since
$\partial M(t_0)/\partial \alpha = \sqrt{2} \neq 0$; this case leads to quadratic heteroclinic tangencies.

\subsection{On the Pikovsky-Topaj  model \cite{PT02} of coupled rotators}
\label{subse:PikTopaj}
%

Let us consider the following system
\begin{equation}
\begin{array}{l}
\dot \psi_1 = 1 - 2 \varepsilon \sin \psi_1 + \varepsilon \sin \psi_2\\
\dot \psi_2 = 1 - 2 \varepsilon \sin \psi_2 + \varepsilon \sin \psi_1
+ \varepsilon \sin \psi_3\\
\dot \psi_3 = 1 - 2 \varepsilon \sin \psi_3 + \varepsilon \sin \psi_2,
\end{array}
\label{Pik1}
\end{equation}
where $\psi_i\in [0,2\pi), i=1,2,3$, are cyclic variables. Thus, the phase space of (\ref{Pik1}) is
the $3$-dimensional torus $\mathbb{T}^3$. System (\ref{Pik1}) is reversible with
respect to the involution ${\cal R}$:
$\; \psi_1 \to \pi- \psi_3 \;\; , \;\; \psi_2 \to \pi - \psi_2 \;\; , \;\; \psi_3 \to \pi -
\psi_1.$
%

System (\ref{Pik1}) was suggested by Pikovsky and Topaj in the paper~\cite{PT02} as a simple model describing the dynamics of $4$ coupled elementary rotators. By means of the coordinate change
\[
\xi=\frac{\displaystyle \psi_1 - \psi_3}{\displaystyle 2}, \;\;  \eta=\frac{\displaystyle \psi_1 +
\psi_3 - \pi}{\displaystyle 2}, \;\; \rho=\frac{\displaystyle \psi_1 + \psi_3 - \pi}{\displaystyle
2}+\psi_2-\pi
\]
and the change in time  $d\tau=dt(2+\varepsilon \cos(\rho-\eta))$ system (\ref{Pik1}) is led into
\begin{eqnarray}
\label{Pik4} \dot \xi &=& \frac{\displaystyle 2 \varepsilon \sin \xi \sin \eta}
{\displaystyle2 + \varepsilon \cos (\rho -\eta)} \nonumber \\
\dot \eta &=& \frac{\displaystyle 1 - \varepsilon \cos(\rho-\eta)  - 2 \varepsilon \cos \xi \cos
\eta}
{\displaystyle 2 + \varepsilon \cos (\rho -\eta)}  \\
\dot \rho &=& 1 \nonumber
\end{eqnarray}
Then time-$1$ Poincar\'e map of system~(\ref{Pik4}) is also
reversible with respect to the same involution $R:\;\;\xi\to\xi,\;\eta\to -\eta$.

It was found in \cite{PT02} that, for small $\varepsilon$, system (\ref{Pik1}) behaves itself as a conservative system close to integrable one and several invariant curves could be observed.
However, when one increases the value of $\varepsilon$ invariant curves break down
and chaos appears (which is already noticed, for instance, at $\varepsilon \approx 0.3$). This picture looks to be quite similar to the conservative case. However, certain principal differences take place. In particular, a ``strange behaviour'' of the invariant measure is observed. Iterations of the initial measure
are convergent to some suitable limit. However,
the limits $t\to +\infty$ and $t\to -\infty$ for the same initial measure are different ( numerically observed, for instance, for values of $\varepsilon \approx 0.3$). This situation is impossible when the invariant measure is absolutely continuous. Therefore, it must be singular and concentrated on "attractors and conservators" at $t\to +\infty$ or "on repellers and conservators" at $t\to -\infty$. Here under the term ``conservator'' we mean the set of self-symmetric non-wandering orbits. Moreover, $+\infty$- and $-\infty$-invariant measures look like
symmetric (with respect to the fixed line of the involution) and having non-empty intersection so there are no gaps between asymmetric and symmetric parts. This means that ``visually'' attractors and repellers intersect and it is an evidence of mixed dynamics in this model.

Moreover, a transition from conservative dynamics to non-conservative one can be generated by bifurcations of periodic orbits. For small enough $\varepsilon$ periods of all such orbits are large and
the corresponding resonance zones are narrow. When increasing $\varepsilon$, periodic orbits of no too large period appear and dissipative phenomena can become observable.
For example, the map $T$ under consideration has no points of period $1$ and $2$ for $\varepsilon<0.6$ but it has, at $\varepsilon=\varepsilon^* \approx  0.445$,
two period $3$ orbits. Notice that these orbits are different since map $T$ has the symmetry $\xi \to 2\pi-\xi$ that implies the appearance of $2$ (in fact, an even number) different orbits. Thus, the scenario is the following:
there is no fixed point for $T^3$ at $\varepsilon<\varepsilon^*$; at
$\varepsilon = \varepsilon^*$ two fixed points with double multiplier $+1$ appear in
$\mbox{Fix $\,R$}$ one symmetric of each other;
at $\varepsilon>\varepsilon^*$ all these orbits fall into four orbits, two symmetric elliptic and two asymmetric saddle. Moreover, the latter orbits satisfy that the Jacobian is greater than $1$ at one point and less than $1$ at other point.

Bifurcations of such type (i.e., having a single point falling into $4$ points) are not typical in one-parameter families even in the reversible case. Here, general bifurcations are met (for symmetric fixed points) of types ``$0\to 2$'' or ``$1\to 3$'', that is, ``conservative'' fold and ``reversible'' pitchfork, respectively. The presence of a typical bifurcation ``$0\to 4$'' says us about the existence of a certain additional degeneracy in the system. The ``clear symmetry'' $\xi\to -\xi$ is not suitable for this r\^ole.
However, system (\ref{Pik4}) possesses such a ``hidden symmetry'' which implies that the map
$T_{(\rho=0) \to (\rho=2 \pi)}$ is the second power of some non-orientable map. This peculiarity is
caused by the fact that the maps
$T_{(\rho=\pi) \to (\rho=2 \pi)}$ and $T_{(\rho=0) \to (\rho=\pi)}$ are conjugate. In particular, one can check that
\begin{equation}
T_{(\rho=\pi) \to (\rho=2 \pi)} = S^{-1} \,  T_{(\rho=0) \to (\rho=\pi)} \, S, \label{sd0vig}
\end{equation}
through the linear change of coordinates
$\xi \to \pi - \xi,$ $\eta \to \eta + \pi$, $\rho \to \rho + \pi$.
Indeed, after this coordinate transformation, the right sides of system (\ref{Pik4}) remain the same, but the limits of integration (along orbits of system (\ref{Pik4}) to get the correspondence map between sections $\rho=a$ and $\rho=b$) are shifted in $\pi$.
Such a property is called {\em time-shift symmetry}.

>From (\ref{sd0vig}) it follows that
$T_{(\rho=0) \to (\rho=2 \pi)}=T_{(\rho=0) \to (\rho=\pi)} \, S \, T_{(\rho=0) \to (\rho=\pi)} \, S^{-1}$. Since $S^2=\mbox{Id}$, one has that $S=S^{-1}$ and, therefore,
\begin{equation}
\begin{array}{l}
T_{(\rho=0) \to (\rho=2 \pi)}=(T_{(\rho=0) \to (\rho=\pi)} S)^2
\end{array}
\label{quadr}
\end{equation}
This means that the map $T_{(\rho=0) \to (\rho=2 \pi)}$ considered is the second power of some map.
Notice that the transformation associated to $S$ is non-orientable and, thus, the map $T_{(\rho=0) \to (\rho=\pi)} S$ is non-orientable as well and, on its turn, our
first-return map $T$ is also the second power of some non-orientable map.

It is straightforward to check that the map $T_{(\rho=0) \to (\rho=\pi)}$ is reversible with respect
to the involution $R_1(x,y)=(-x,-y)$ and that the map $T_{(\rho=0) \to (\rho=\pi)} \,S$ is reversible under the involution $R(x,y)=(x,-y)$.
Thus, the bifurcation of map $T^3$ at $\varepsilon=\varepsilon^*$ can be treated  as a bifurcation of a fixed point with multipliers $(+1,-1)$ in the case of a non-orientable map (in fact, the map $(T_{(\rho=0) \to (\rho=\pi)}\, S)^3$).
So, summarising, in our case this bifurcation leads to the appearance of two elliptic points of period $2$ on $\mbox{Fix$\,R$}$ and a symmetric couple of saddle fixed points
(that is, outside $\mbox{Fix $\,R$}$ and symmetric one to each other).
These saddle fixed points are not conservative. It can be
checked numerically that the Jacobian of one point is greater than $1$ and less than $1$ at other point. Due to reversibility, the stable and unstable manifolds of saddles pairwise intersect and form a ``heteroclinic tangle'' zone. This zone is extremely narrow since the separatrix splitting is exponentially small. However, moving slightly away from the bifurcation moment we can find numerically heteroclinic tangencies and, hence, moments of creation of non-transversal heteroclinic cycles.
Since the saddles involved are not conservative, it follows from~\cite{LSt} the phenomenon of mixed dynamics.

\section{Cross-form type equations for reversible maps. Proof of Lemmas~\ref{LmSaddle}
and~\ref{LmLocalMap}}
\label{se:crossL1L2}

\subsection{Cross-form for reversible maps}
\label{subse:cross_form}

As it will be seen along this section, the so-called Shilnikov \emph{cross-form} variables
constitute an essential (and natural) tool to deal with reversible maps and a simple way to generate them. The first part will be devoted to introduce such variables and to present some of its main characteristics. In the second part we apply them to prove Lemmas~\ref{LmSaddle}
and~\ref{LmLocalMap}.

We say that a map is in \emph{cross-form} if it is written as
\[
\left\{
\begin{array}{l}
\bar{x} = h(x,\bar{y}), \\
     y  = h(x,\bar{y}),
\end{array}
\right.
\]
On the other hand, let us
consider a diffeomorphism $F$ of the plane which is reversible with respect to a (in general, non-linear) involution $\Rrev$ ($\Rrev^2=\mbox{id}$, $\Rrev\neq \mbox{id}$), having dim Fix$\Rrev=1$. Moreover, let assume that the involution $\Rrev$ reverses orientation (that is, $\mbox{det$\,D\Rrev$}<0$), which is the most common situation in the literature.
Our aim is to show how reversible maps can be expressed in cross-form type equations
and, conversely, how reversibility can be derived from this form.

As an starting point, let us consider the linear set up, that is, when the reversor $\Rrev$ is the linear involution $\Lrev:(x,y) \mapsto (y,x)$. In this case, the following result holds:
\begin{lm}
\label{appendix:lemma:1} Any diffeomorphism $F:(x,y)\mapsto (\bar{x},\bar{y})$ defined, implicitly, by means of equations of type
\begin{equation}
\label{appendix:F_rev} F: \ \left\{
\begin{array}{l}
\bar{x} = f(x,\bar{y}), \\
     y  = f(\bar{y},x)
\end{array}
\right.
\end{equation}
is always reversible with respect to $L(x,y)=(y,x)$.
\end{lm}

\proof Remind that if $G$ is a ${\Lrev}$-reversible diffeomorphism it must satisfy that $G \circ
{\Lrev} \circ G = {\Lrev}$ or, equivalently, ${\Lrev} \circ G \circ {\Lrev} = G^{-1}$ or $\left(
{\Lrev} \circ G \circ {\Lrev} \right)^{-1} = G$. In our case we will prove that $F$, defined
by~(\ref{appendix:F_rev}), verifies the latter relation for $G=F$ and, consequently, is
$\Lrev$-reversible. To do it, we will use an equivalent expression for the inverse of a planar
diffeomorphism. Precisely, if
\[
H: \ \left\{
\begin{array}{l}
 \bar{x} = h_1(x,y), \\
 \bar{y} = h_2(x,y),
\end{array}
\right.
\]
the corresponding inverse map $H^{-1}: (x,y) \mapsto (\bar{x},\bar{y})$ can be implicitly written
through the expression
\[
H^{-1}: \ \left\{
\begin{array}{l}
     x = h_1(\bar{x},\bar{y}), \\
     y = h_2(\bar{x},\bar{y}).
\end{array}
\right.
\]
This is clear since $(x,y)=H(\bar{x},\bar{y})$ implies that $(\bar{x},\bar{y})=H^{-1}(x,y)$. An
\emph{algorithmic} way to get it consists in swapping bars among the variables, that is $x
\leftrightarrow \bar{x}$ and $y \leftrightarrow \bar{y}$. We apply this procedure to compute formally an expression for $( {\Lrev} \circ F \circ {\Lrev})^{-1}$ and to check
afterwards that it coincides with $F$. Let us compute it step by step. First we have
\[
F\circ {\Lrev}: \ \left\{
\begin{array}{l}
\bar{x} = f(y,\bar{y}), \\
     x  = f(\bar{y},y).
\end{array}
\right.
\]
To apply $\Lrev$ onto $F \circ {\Lrev}$ corresponds to swap $\bar{x} \leftrightarrow \bar{y}$ in
the precedent expression:
\[
{\Lrev} \circ F\circ {\Lrev}: \ \left\{
\begin{array}{l}
\bar{y} = f(y,\bar{x}), \\
     x  = f(\bar{x},y).
\end{array}
\right.
\]
And finally, to get its inverse we swap bars and no-bars, that is $x \leftrightarrow \bar{x}$ and
$y \leftrightarrow \bar{y}$. Performing this change we obtain
\[
\left({\Lrev} \circ F\circ {\Lrev}\right)^{-1}: \ \left\{
\begin{array}{l}
     y  = f(\bar{y},x), \\
     \bar{x}  = f(x,\bar{y}),
\end{array}
\right.
\]
which is exactly the expression for $F$. So, $F$ given in the form~(\ref{appendix:F_rev}) is always ${\Lrev}$-reversible.

\qed

\noindent This result can be useful to provide suitable local expressions for reversible
diffeomorphisms in the plane. Thus we have the following lemma.
\begin{lm}
Let $F=(f_1,f_2)$ be a planar diffeomorphism, reversible with respect a general
involution $\Rrev$,  ${\cal C}^{r}$, $r\geq 1$,
orientation reversing and with $\mbox{dim\,Fix$\,\Rrev$}=1$.
Let us assume the origin $(0,0)$ a fixed point of the involution $\Rrev$, that is $(0,0) \in \mbox{Fix$\Rrev$}$.

Then, if $D_{xx}f_1  + D_{yy} f_2 \ne 0$ at $(0,0)$ there exist local coordinates, that we denote again by $(x,y)$, in which $F$ admits the following implicit (normal) form
\[
\left\{
\begin{array}{rcl}
\bar{x} &=& g(x,\bar{y}), \\
     y &=& g(\bar{y},x).
\end{array}.
\right.
\]
This map is reversible with respect to $L(x,y)=(y,x)$.
\end{lm}
In the case of a saddle fixed point, the concrete type of implicit normal form that can be obtained is given in equation~(\ref{CrossFormSaddle}).

\proof This will be achieved in two steps:
\begin{itemize}
\item[$(i)$] First we apply Bochner Theorem~\cite{B45} which allows us to conjugate, around $(0,0)$, our involution $\Rrev$ to its linear part $DG|_{(0,0)}$.
\item[$(ii)$] Using that the partial derivatives on $(0,0)$ do not vanish simultaneously, we apply Implicit Function Theorem to reach the final form.
\end{itemize}
We proceed as follows:

\begin{itemize}
\item[$(i)$] Notice that if $\Rrev$ is a (general) involution and $p \in \mbox{Fix
$\,\Rrev$}$ then its linear part $D\Rrev|_p$ is an involution as well. Indeed,
\[
\mbox{Id}=\Rrev^2 \Rightarrow I = D(\Rrev^2)(p)=D\Rrev|_{\Rrev(p)} \cdot D\Rrev|_{p} = \left(
D\Rrev|_{p} \right)^2.
\]
Bochner Theorem ensures the existence of a ${\cal C}^{r}$-diffeo $\psi$ which conjugates,
locally around $p$, $\Rrev$ to $D\Rrev|_p$. We include, for completeness, a simple proof
of this fact given in~\cite{RQ92}. From the equality
\begin{eqnarray*}
\lefteqn{D\Rrev|_p \circ \left( \Rrev + D\Rrev|_p \right) = D\Rrev|_p \circ \Rrev + \mbox{id} =} \\
&& \qquad \mbox{id} + D\Rrev|_p \circ G = \left( G + DG|_p \right) \circ G
\end{eqnarray*}
it follows that $DG|_p \circ \left( \Rrev + D\Rrev|_p \right) = \left( \Rrev + D\Rrev|_p \right)
\circ \Rrev$. We define $\psi = \Rrev + D\Rrev|_p \in {\cal C}^r$ and check that it is a diffeomorphism in a neighbourhood of $p$:
\[
D\psi|_p = D\left( \Rrev + D\Rrev|_p \right)|_p = D\Rrev|_p + D\Rrev|_p = 2 D\Rrev|_p
\]
and so $\det \,D\psi|_p = 2 \det \,D\Rrev|_p \neq 0$ since $\Rrev$ is a diffeomorphism around
$p$. So $\psi$ is a ${\cal C}^{r}$-diffeomorphism which conjugates $\Rrev$ to $D\Rrev|_p$ around
$p$.

Since $\Rrev$ is orientation reversing its linear part around $p$, $D\Rrev|_p$, is also orientation
reversing. Following~\cite{RQ92} for instance, we know that there exists a transformation which conjugates $D\Rrev|_p$ to the linear involution $\Lrev(x,y)=(y,x)$, which will be the one we will consider, locally, from now on.

\item[$(ii)$] Let us assume, for instance, that $D_{yy}f \neq 0$ at $(0,0)$.
Using Implicit Function Theorem, we can write from equation $\bar{y} = f_2(x,y)$ an expression for
$y$, say $y = g(\bar{y},x)$, for a suitable function $g$. Substituting it into the equations defining $F$ we get a (locally) equivalent expression for $F$:
\[
F: \, \left\{
\begin{array}{l}
\bar{x} = f_1(x,y) = f_1(x, g(\bar{y},x)) =: h(x,\bar{y}),\\
y = g(\bar{y},x).
\end{array}
\right.
\]
As stated above, we can assume $F$ to be locally conjugated around the origin to the linear
involution $\Lrev:(x,y) \mapsto (y,x)$. So in that variables (to simplify the notation we keep the same name for the variables and the functions involved)  it must satisfy that $({\Lrev} \circ F \circ {\Lrev})^{-1} =
F$. Applying the procedure introduced in Lemma~\ref{appendix:lemma:1}, one obtains that
\[
F \circ {\Lrev} : \, \left\{
\begin{array}{l}
\bar{x} = h(y,\bar{y}), \\
     x = g(\bar{y},y).
\end{array}
\right.
\]
We apply $L$ (that corresponds to swapping $\bar{x}$ and $\bar{y}$,
\[
{\Lrev} \circ F \circ {\Lrev} : \, \left\{
\begin{array}{l}
\bar{y} = h(y,\bar{x}), \\
     x = g(\bar{x},y)
\end{array}
\right.
\]
and, finally, we swap $(x,y)$ for $(\bar{x},\bar{y})$,
\[
\left({\Lrev} \circ F \circ {\Lrev} \right)^{-1} : \, \left\{
\begin{array}{l}
      y = h(\bar{y},x), \\
\bar{x} = g(x,\bar{y}).
\end{array}
\right.
\]
Since it must coincide with $F$ it turns out that $h(x,\bar{y})=g(x,\bar{y})$ and so
\[
F: \, \left\{
\begin{array}{l}
\bar{x} = g(x,\bar{y}), \\
     y = g(\bar{y},x).
\end{array}
\right.
\]
\end{itemize}

\qed
We present now a counterpart result when the map is given in implicit form.
\begin{lm}
\label{appendix:B:implicitrev} Any map $G: (x,y) \mapsto (\bar{x},\bar{y})$ given by
\[
\left\{
\begin{array}{l}
g(x,y,\bar{x},\bar{y}) = 0, \\
g(\bar{y},\bar{x},y,x) = 0,
\end{array}
\right.
\]
is ${\Lrev}$-reversible, where ${\Lrev}:(x,y)\mapsto (y,x)$.
The second equation $g(\bar{y},\bar{x},y,x) = 0$ is a kind of ${\Lrev}$-conjugate
of the first equation $g(x,y,\bar{x},\bar{y})=0$.
\end{lm}

\proof It is enough to check that ${\Lrev} \circ G \circ {\Lrev} = G^{-1}$. To do it we proceed again as in Lemma~\ref{appendix:lemma:1}. First, remind that an implicit expression for $G^{-1}$ is always
obtained by swapping \emph{bars} for \emph{no-bars}, that is, $(x,y)\leftrightarrow
(\bar{x},\bar{y})$. So
\[
G^{-1}: \ \left\{
\begin{array}{l}
g(\bar{x},\bar{y},x,y)=0, \\
g(y,x,\bar{y},\bar{x})=0.
\end{array}
\right.
\]
On the other hand we compute ${\Lrev} \circ G \circ {\Lrev}$. Thus,
\[
G \circ {\Lrev}: \ \left\{
\begin{array}{l}
g(y,x,\bar{x},\bar{y})=0, \\
g(\bar{y},\bar{x},x,y)=0,
\end{array}
\right.
\]
and, swapping $(\bar{x},\bar{y})$ for $(\bar{y},\bar{x})$, we get
\[
{\Lrev} \circ G \circ {\Lrev}: \ \left\{
\begin{array}{l}
g(y,x,\bar{y},\bar{y})=0, \\
g(\bar{x},\bar{y},x,y)=0,
\end{array}
\right.
\]
which coincides with $G^{-1}$. Therefore the lemma is proved.

\qed

The following result establishes an interesting relation between polynomial
reversible and area preserving maps.
\begin{lm}[\cite{R97}]
\label{appendixB:polynrevr:conserv}
Any Taylor truncation of a planar polynomial diffeomeorphism which is reversible with respect to a linear involution is area preserving. In particular, this applies to the truncation of a normal form of such diffeomorphisms.
\end{lm}

\proof Let $\bar{z}=G(z)$ a polynomial planar map which is reversible with respect to a linear
involution $S$ ($S^2=\Id$, $S\neq \Id$). This means that $S \circ G \circ S = G^{-1}$ and, in
particular, that $G^{-1}$ is also a polynomial. Differentiating the latter expression we get
\begin{eqnarray*}
&& S \, DG|_{Sz} \, S = D(G^{-1})|_z= \left( DG|_{G^{-1}(z)} \right)^{-1} \Rightarrow \\
&& \det (S \, DG|_{Sz} \, S) = \frac{1}{\det DG|_{G^{-1}(z)}}.
\end{eqnarray*}
Using that $\det (S\, DG|_{Sz} \, S)=(\det S)^2 \det DG|_{Sz} = \det DG|_{Sz}$ it follows
that
\begin{equation}
\label{appendixB:polynrevr:conserv:proof} (\det DG|_{Sz}) \cdot (\det DG|_{G^{-1}(z)} ) = 1, \quad
\forall z.
\end{equation}
Since $G$ and $G^{-1}$ are polynomials and $S$ linear we obtain that $\det DG|_{Sz}$ and $\det
DG|_{G^{-1}(z)}$ are polynomials as well. But the product of two polynomials is a constant if and
only if they are constant, that is, $\det DG|_z \equiv k=\mbox{constant}$. Thus,
from~(\ref{appendixB:polynrevr:conserv:proof}) it follows that $k^2=1$ and, therefore, $\det
DG|_z=\pm 1$, $\forall z$.

\qed

And last, but not least, we remark another interesting property regarding this cross-form type: any polynomial truncation of a reversible diffeomorphism written in cross-form type is also in cross-form type and, consequently, it is reversible. This means, from Lemma~\ref{appendixB:polynrevr:conserv}, that this truncation is always area-preserving.

\subsection{Proof of Lemma~\ref{LmSaddle}}
\label{PrL1}

Let $O$ be a fixed saddle
point of a reversible map $T_0$. Applying Bochner Theorem~\cite{MZ55}, we can assume the existence of local coordinates around $O$ such that $O$ is located at the origin and that the involution $R$ is exactly $(x,y) \mapsto (y,x)$ in these coordinates.

Let $x=\nu(y)$ be the equation of the stable manifold. Then, by
the $R$-reversibility, $y=\nu(x)$ is the equation of the unstable
manifold. If $|d\nu/dy|<1$, we perform the transformation
$x_{new} = x - \nu(y)$, $y_{new} = y - \nu(x)$,
while, if $|d\nu/dy|>1$, the change is
$x_{new} = y - \nu(x)$, $y_{new} = x - \nu(y)$.
After such transformation, which commutes with $R$, the equations of the stable and unstable manifolds become $y=0$ and $x=0$, respectively. Thus, in the corresponding  local coordinates, the map can be represented in the following form
\begin{equation}
\bar x = \lambda x  + g_1(x,y),\qquad \bar y = \lambda^{-1} y + g_2(x,y)
\label{SF1}
\end{equation}
where $g_1(0,y)\equiv 0,g_2(x,0)\equiv 0$  and
$g_i^\prime(0,0)=0,i=1,2$. It is very convenient to rewrite this
equation in the so-called {\em cross-form}:
\begin{equation}
\begin{array}{l}
\bar x = \lambda x  + \tilde g_1(x,\bar y),\qquad y = \lambda \bar y + \tilde g_2(x,\bar y)
\label{SF2}
\end{array}
\end{equation}
Equation~(\ref{SF2}) comes from (\ref{SF1}) writing $y = F(x,\bar y)$ (which exists due to the Implicit Function Theorem) and substituting it into the first equation:
$\bar x = \lambda x  + g_1(x,F(x,\bar y))$.
The $R$-reversibility of (\ref{SF2}) implies that $\tilde g_1(x,y) \equiv \tilde g_2(y,x)$ so we can represent map (\ref{SF2}) in the form
\begin{equation}
\begin{array}{l}
\bar x = \lambda x  + \varphi_1(x)+ \psi_1(\bar y)x +  \rho_1(x,\bar y) x^2 \bar y\\
y = \lambda \bar y + \varphi_1(\bar y)+ \psi_1(x)\bar y + \rho_1(\bar y,x)x \bar y^2 \label{SF3}
\end{array}
\end{equation}
Perfoming the $R$-invariant change of variables
\begin{equation}
\xi = x + x h_1(y),\qquad \eta = y + y h_1(x)
\label{RInSaCh1}
\end{equation}
with $h_1(0)=0$, it turns out the following equation for $\bar\xi$:
\begin{eqnarray*}
\bar\xi &=& \bar x + \bar x h_1(\bar y) =  \nonumber \\
&&
\lambda\xi -  xh_1(y) + x\psi_1(\bar y ) +
(\lambda x + \psi_1(\bar y)x + \varphi_1(x))h_1(\bar y) + \\
&& \varphi_1(\xi)+  O(\xi^2\bar\eta) = \\
&& \lambda\xi +   \varphi_1(\xi)+ O(\xi^2\bar y) + \nonumber \\
&& x\left[-h_1(\lambda \bar y + \varphi_1(\bar y)) + \psi_1(\bar y) +
(\lambda + \psi_1(\bar y))h_1(\bar y)\right].\nonumber
\label{FuEq1}
\end{eqnarray*}
Since we want the expression in the square brackets to vanish identically, we ask
the function $h_1(y)$ to satisfy the functional equation
\begin{equation}
\label{h1eq}
  h_1(\lambda\bar y + \varphi_1(\bar y))=h_1(\bar y)(1+\lambda^{-1}\psi_1(\bar y)) +\lambda^{-1}\psi_1(\bar y),
\end{equation}
which has solutions $h_1=h_1(u)$ in the class of $C^{r-1}$-functions.
Indeed, we can consider (\ref{h1eq}) as an equation
for the strong stable invariant manifold of the following planar
map
\begin{equation*}
\label{eqh1man}
\begin{array}{rcl}
\bar h_1 &=& h_1(1+\lambda^{-1}\psi_1(u)) +\lambda^{-1}\psi_1(u), \\
\bar u   &=& \lambda u + \varphi_1(u).
\end{array}
\end{equation*}
Since $0<|\lambda|<1$, $\psi_1(0)=0$ and $\varphi_1(0)=\varphi_1^\prime
(0)=0 $,  this map has strong stable invariant manifold
$W^{ss}$ passing through the origin, that is, satisfying an equation $h_1=h_1(u)$ with $h_1(0)=0$.
Therefore, after the $R$-invariant change (\ref{RInSaCh1}), the map (\ref{SF3}) takes the form
\begin{equation}
\begin{array}{rcl}
\bar x &=& \lambda x  + \varphi_1(x)+  \rho_2(x,\bar y) x^2 \bar y\\
y &=& \lambda \bar y + \varphi_1(\bar y) + \rho_2(\bar y,x)x \bar y^2
\label{SF4}
\end{array}
\end{equation}
Applying a $R$-invariant change of variables of the form
\begin{equation*}
\xi = x + h_2(x)x, \qquad \eta = y + h_2(y)y
\label{RInSaCh3}
\end{equation*}
with $h_2(0)=0$, the first equation of system (\ref{SF4}) can be
rewritten, in these new coordinates, as follows
\begin{eqnarray}
\bar\xi &=& \lambda\xi + x\left[-\lambda h_2(x)
+\tilde\varphi_1(x) + \right. \label{SFs5}  \\
&&
\left. h_2(\lambda x + \varphi(x))(\lambda +
\tilde\varphi_1(x))\right] + O(\xi^2\eta), \nonumber
\end{eqnarray}
where we have denoted $\varphi_1(x)\equiv \tilde\varphi_1(x) x$. As we did above for $h_1$, we seek for a function $h_2$ satisfying the following equation
\begin{equation}
h_2(\lambda x + \varphi(x))= (1 +\lambda^{-1}\tilde\varphi_1(x))^{-1}(h_2(x) -
\lambda^{-1}\tilde\varphi_1(x)),
\label{h2n}
\end{equation}
which vanishes the expression inside the square brackets
in (\ref{SFs5}). As before, equation (\ref{h2n}) has solutions $h_2=h_2(u)$ in the class
of $C^{r-1}$-functions. Again, one can consider the
expression (\ref{h2n}) as an equation for the strong stable
invariant manifold associated to the following planar map
\begin{equation*}
\label{eqh2man}
\begin{array}{l}
\bar h_2 =  (1 +\lambda^{-1}\tilde\varphi_1(x))^{-1}(h_2 -
\lambda^{-1}\tilde\varphi_1(x)),  \\
\bar u = \lambda u + \varphi_1(u).
\end{array}
\end{equation*}
Having in mind that $0<|\lambda|<1$ and $\varphi_1(0)=\varphi_1^\prime (0)=0$,
this map admits strong stable invariant manifold $W^{ss}$ passing
through the origin, i.e., having an equation $h_2=h_2(u)$ with
$h_2(0)=0$. This completes the proof of the Lemma.

\subsection{Proof of Lemma~\ref{LmLocalMap}}
\label{subse:proof_lemma_LocalMap}
We write the map $T_0$ in the following form
\begin{equation*}
\label{v0}
\begin{array}{l}
\bar x \; = \; \lambda x + \hat{h}(x,y), \;\;\; \bar y \; = \; \gamma y + \hat{g}(x,y)
\end{array}
\end{equation*}
where we assume that $\gamma = \lambda^{-1}$ and
\[
\hat{h}(x,y) \equiv x^2y(\beta_1 + O(|x|+|y|)), \quad
\hat{g}(x,y) \equiv xy^2(\beta_2 + O(|x|+|y|)).
\]
Consider the following operator $\Phi:[(x_j,y_j)]^k_{j=0} \mapsto [(\bar x_j, \bar y_j)]^k_{j=0}$:
\begin{equation}
\label{v1}
\begin{array}{l}
\bar x_j =\lambda^j x_0 + \sum\limits^{j-1}_{s=0}\lambda^{j-s-1}\hat{h}(x_s,y_s,\mu), \\
\bar y_{j} = \gamma^{j-k}y_k - \sum\limits^{k-1}_{s=j}\gamma^{j-s-1}\hat{g}(x_s,y_s,\mu),
\end{array}
\end{equation}
where
$j=0,1,\dots,k$. The operator $\Phi$ is defined on the set
\[
Z(\delta) \; = \; \{z=[(x_j,y_j)]^k_{j=0}, \; \Vert z \Vert \le \delta \} \;\;,
\]
where the norm $\Vert \cdot \Vert $ is given as the maximum of modulus of components  $x_j,y_j$ of
the vector $z$. Notice that if $z_0 \; = \; [(x^0_j,y^0_j)]^k_{j=0}$ is a fixed point of $\Phi$, then
the following diagram takes place
\[
(x^0_0,y^0_0) \stackrel{T_0}{\longrightarrow}(x^0_1,y^0_1) \stackrel{T_0}{\longrightarrow}\dots
\stackrel{T_0}{\longrightarrow}(x^0_k,y^0_k),
\]
i.e. the fixed point of $\Phi$ gives a segment of an orbit of $T_0$.

It is known \cite{AS73} that, for small enough $\delta
=\delta _0$ and $\vert x_0\vert \le \delta _0/2,$ $\vert y_k\vert
\le \delta _0/2$, the operator $\Phi$ maps the set $Z(\delta _0)$
into itself and  it is contracting. Thus, map (\ref{v1}) has a
unique fixed point  $z_0 \; = \;
[(x^0_j(x_0,y_k),y^0_j(x_0,y_k)]^k_{j=0}$ that is limit of
iterations under $\Phi$ for any initial point from $Z(\delta_0)$.
Thus, the coordinates $x^0_j$ and $y^0_j$ can be found by applying successive approximations. As an initial approximation, we
take the solution of the linear problem:
\[
x^{0(1)}_j = \lambda^j x_0, \qquad  y^{0(1)}_j = \gamma^{j-k}y_k
\]
It follows from (\ref{v1}) that the second approximation has a form
\[
\begin{array}{l}
x^{0(2)}_j = \lambda^j x_0 + \sum\limits^{j-1}_{s=0}\lambda^{j-s-1} \lambda^{2s}\gamma^{s-k} x_0^2 y_k \times \\
\qquad (\beta_1 + O(|\lambda|^s |x_0| + |\gamma|^{s-k} |y_k|)) =  \\
\; \lambda^j x_0 + \lambda^{j}\gamma^{-k}\sum\limits^{j-1}_{s=0}\lambda^{-1} \lambda^{s}\gamma^{s}
x_0^2 y_k (\beta_1 + O(|\lambda|^s |x_0| + |\gamma|^{s-k} |y_k|)) = \\
\; \lambda^j x_0 + (j-1)\lambda^{j}\gamma^{-k}\lambda^{-1} x_0^2 y_k \left(\beta_1 +
O(|\lambda|^s |x_0| + |\gamma|^{s-k} |y_k|)\right) ,\\
\\
y^{0(2)}_j = \gamma^{j-k}y_k + \sum\limits^{k-1}_{s=j}\gamma^{j-s-1}
\lambda^{s}\gamma^{2(s-k)} x_0 y_k^2 \times \\
\qquad (\beta_2 + O(|\lambda|^s |x_0| + |\gamma|^{s-k} |y_k|)) = \\
\;\gamma^{j-k}y_k + \gamma^{j-2k} \sum\limits^{k-1}_{s=j}\gamma^{-1} \lambda^{s}\gamma^{s} x_0
y_k^2 (\beta_2 + O(|\lambda|^s |x_0| + |\gamma|^{s-k} |y_k|)) = \\
\;\gamma^{j-k}y_k +
(k-j)\gamma^{j-2k-1} x_0 y_k^2 (\beta_2 + O(|\lambda|^s |x_0| + |\gamma|^{s-k} |y_k|))
\end{array}
\]
Since $\gamma =\lambda^{-1}$, it follows from
the precedent expression that
\begin{equation}
\label{v02}
\begin{array}{l}
|x_j^{0(2)} - \lambda^j x_0| \leq L_1 j\lambda^{j+k}, \\
|y_j^{0(2)} - \lambda^{k-j}y_k|  \leq L_2 (k-j)\lambda^{2k-j},
\end{array}
\end{equation}
where $L_1$ and  $L_2$ are some positive constant independent of  $j$ and $k$. Substituting (\ref{v02}) into (\ref{v1}) as the initial approximation, then the following ones will also satisfy estimates (\ref{v02}), with the same constants $L_1$ and $L_2$. Thus, formula
(\ref{LocalMapk}) is valid for the coordinates $x_l^0$ and $y_0^0$, fixed point of $\Phi$.

The estimates for the derivatives of the functions $x_l^0$ and $y_0^0$
are deduced in the same way as done in~\cite{GS92} (see also modified versions of the
proof in~\cite{book,GST07,GST08}).

\section{Proofs of Lemmas~\ref{lm:frm},~\ref{lmnsym} and~\ref{lmgenel}.}
\label{se:technical_lemmas}

\subsection{Proof of Lemma~\ref{lm:frm}}

Since coordinates $(x_{01},y_{01})$ on $\sigma_k^{01}$ are uniquely
determined via cross-coordinates $(x_{01},y_{11})$ in equations~(\ref{T1k}), we can express $T_{km}$ as a map defined on points $(x_{01},y_{11})$ and acting by the rule  $(x_{01},y_{11})\mapsto (\bar
x_{01},\bar y_{11})$. As a result of this, we can express the map $T_{km}$ in the following form
\[
\begin{array}{l}
x_{02} - x_2^+ = a \lambda_1^k x_{01} + b (y_{11}-y_1^-) + l_{02}(y_{11}-y_1^-)^2 + \\ \qquad \tilde\varphi_{1k}(x_{01},y_{11}-y_1^-,\mu),\\[1.2ex]
\lambda_2^m y_{12}\left(1+ m\lambda^m_2h_m^2(y_{12},x_{02},\mu)\right) = \mu + c \lambda_1^k x_{01}
+ d (y_{11}-y_1^-)^2 + \\
\qquad f_{11} \lambda_1^k
x_{01}(y_{11}-y_1^-)+f_{03}(y_{11}-y_1^-)^3+ \tilde\varphi_{2k}(x_{01},y_{11}-y_1^-,\mu), \\[1.2ex]
\lambda_2^m x_{02}\left(1+ m\lambda^m_2h_m^2(x_{02},y_{12},\mu)\right)  = \mu + c \lambda_1^k \bar y_{11} + d (\bar x_{01}-y_1^-)^2 + \\
\qquad f_{11} \lambda_1^k \bar y_{11}(\bar x_{01}-y_1^-) + f_{03}(\bar x_{01}-y_1^-)^3 + \tilde\varphi_{2k}(y_{11}-y_1^-,x_{01},\mu),\\[1.2ex]
y_{12} - x_2^+ = a \lambda_1^k \bar y_{11} + b (\bar x_{01}-y_1^-) + l_{02}(\bar x_{01}-y_1^-)^2 + \\
\qquad \tilde\varphi_{1k}(y_{11}-y_1^-,x_{01},\mu),\\
\end{array}
\]
where the coordinates $x_{02}$ and $y_{12}$ are ``intermediate'' and
\begin{equation}
\begin{array}{l}
\tilde\varphi_{1k}(u,v,\mu) = O\left(\lambda_1^{2k}u^2 +
|\lambda_1^{k}||uv|
+|v|^3\right), \\
\tilde\varphi_{2k}(u,v,\mu) = O\left(\lambda_1^{2k}(u^2 + |u|v^2)
+ |\lambda_1^{k}||u|v^2\right) + o(v^3).
\end{array}
\label{est*}
\end{equation}
Now we perform the following shift in the coordinates
\begin{eqnarray*}
\xi_1 &=& x_{01}-x_1^+ + \nu_{km}^1, \qquad \eta_1=y_{11}-y_1^- +
\nu_{km}^1, \\
\xi_2 &=& x_{02}-x_2^+ + \nu_{km}^2, \qquad
\eta_1=y_{12}-y_2^- + \nu_{km}^2,
\end{eqnarray*}
where $\nu_{km}^i=O(\lambda_1^k)$, $i=1,2,$ are some small
coefficients which does not destroy the reversibility due to the
condition (\ref{xeqy}). Then, for suitable $\nu_{km}^i$,  map $T_{km}$
becomes
\begin{equation}
\fl
\begin{array}{l}
\xi_2 = a \lambda_1^k \xi_1 + b\eta_1 + l_{02}\eta_1^2 + \tilde\varphi_{1k}(\xi_1,\eta_1,\mu),\\[1.2ex]
\lambda_2^m \eta_2\left(1+ m\lambda^m_2h_m^2(\eta_2,\xi_2,\mu)\right) = \\
\qquad  \tilde\mu + c \lambda_1^k \xi_1 + d \eta_1^2  +
f_{11}\lambda_1^k\xi_1\eta_1 + f_{03}\eta_1^3+ \tilde\varphi_{2k}(\xi_1,\eta_1,\mu),
\\[1.2ex]
\lambda_2^m \xi_{2}\left(1+ m\lambda^m_2h_m^2(\xi_2,\eta_2,\mu)\right) = \\
\qquad  \tilde\mu + c \lambda_1^k \bar\eta_1 + d \bar\xi_1^2
f_{11}\lambda_1^k\bar\eta_1\bar\xi_1 + f_{03}\bar\xi_1^3 + \tilde\varphi_{2k}(\bar\eta_1,\bar\xi_1,\mu),\\[1.2ex]
\eta_2 = a \lambda_1^k \bar\eta_1 + b\bar\xi_1 + l_{02}\bar\xi_1^2 +
\tilde\varphi_{1k}(\bar\eta_1,\bar\xi_1,\mu),\\
\end{array}
 \label{Tkm+0}
\end{equation}
where, since relation~(\ref{xeqy}) holds, we have that $\tilde\mu = \mu - \lambda_2^m (\alpha_2^*+\dots) + c \lambda_1^k(\alpha_1^*+\dots)$ and the new
functions $\tilde\varphi_{1k}$ and $\tilde\varphi_{2k}$ satisfy again conditions (\ref{est*}).
One must, however, to consider coefficients $a,b,...,f_{03}$ in (\ref{Tkm+0}) to be shifted by values of order $O(k\lambda_1^k)$ when comparing them
with the initial coefficients in~(\ref{T12}).
Substituting into the second and third equations of (\ref{Tkm+0}) the expressions for $\xi_2$ and $\eta_2$ given by the first and the fourth precedent equations, we get an expression for $T_{km}$ of form
\begin{equation*}
\begin{array}{l}
\lambda_2^m\left(a \lambda_1^k \bar\eta_1 + b\bar\xi_1 + l_{02}\bar\xi_1^2 \right) = \\ \qquad \tilde\mu + c\lambda_1^k \xi_1 + d \eta_1^2  + f_{11}\lambda_1^k\xi_1\eta_1 +
f_{03}\eta_1^3+ \tilde\varphi_{2k}, \\[1.2ex]
\lambda_2^m \left(a\lambda_1^k \xi_1 + b\eta_1 + l_{02}\eta_1^2 \right) = \\
\qquad \tilde\mu + c \lambda_1^k \bar\eta_1 + d \bar\xi_1^2 +
f_{11}\lambda_1^k\bar\eta_1\bar\xi_1 + f_{03}\bar\xi_1^3 + \tilde\varphi_{2k},\\
\end{array}
\label{Tkm+1}
\end{equation*}
which can be rewritten as
\begin{equation}
\label{Tkm+2}
\begin{array}{l}
a\lambda_2^m \bar\eta_1 + b\lambda_2^m\lambda_1^{-k}\bar\xi_1 +
l_{02}\lambda_2^m\lambda_1^{-k}\bar\xi_1^2 = \\
\qquad \tilde\mu \lambda_1^{-k} + c \xi_1 + d\lambda_1^{-k}\eta_1^2  + f_{11}\xi_1\eta_1 +
f_{03}\lambda_1^{-k}\eta_1^3+ \lambda_1^{-k}\tilde\varphi_{2k}, \\[1.2ex]
a\lambda_2^m\xi_1 + b\lambda_2^m\lambda_1^{-k}\eta_1 + l_{02}\lambda_2^m
\lambda_1^{-k}\eta_1^2  =  \\
\qquad \tilde\mu \lambda_1^{-k} + c \bar\eta_1 + d\lambda_1^{-k} \bar\xi_1^2 +
f_{11}\bar\eta_1\bar\xi_1 + f_{03}\lambda_1^{-k}\bar\xi_1^3 + \lambda_1^{-k}\tilde\varphi_{2k}.
\end{array}
\end{equation}
Notice that the functions $\tilde\varphi_{2k}$ here may be changed in comparison with those in (\ref{Tkm+0}) but still fulfill relations (\ref{est*}).
Finally, rescaling coordinates,
\[
\xi_1 = \lambda_1^k x, \qquad \eta_1=\lambda_1^k y,
\]
system (\ref{Tkm+2}) takes the form (\ref{frmSO}) where the coefficients
$c,d,\dots,l_{02}$ are ``original'' ones (i.e., those appearing in formula~(\ref{T12})).

\subsection{Proof of Lemma~\ref{lmnsym}}
The rescaled form (\ref{frmSO}) of the  first-return
map $T_{km}$ is, of course, implicit one and it corresponds to a
formal representation
${\ds T_{km}: \;\; F(\bar x, \bar y) \equiv G(x,y) }$
%
which can be written in the explicit form
${\ds (\bar x, \bar y) \equiv T_{km}(x,y) \equiv F^{-1}G(x,y)}.$
%
Then we can find the Jacobian of $T_{km}$ using the relation
\begin{equation}
\displaystyle D(T_{km})\bigl|_{(x,y)} \equiv D(F^{-1})\bigl|_{G(x,y)} D(G(x,y)),
\label{Jacm}
\end{equation}
where $D(\cdot)$ is the corresponding (differential) Jacobi
matrix. At the fixed point $(x=\xi^*,y=\eta^*)$ of $T_{km}$ we can
rewrite (\ref{Jacm}) as follows
%
\[
D(T_{km})\bigl|_{(\xi^*,\eta^*)} \equiv
\left(DF\bigl|_{(\xi^*,\eta^*)}\right)^{-1}
DG\bigl|_{(\xi^*,\eta^*)}.
\]
%
We find from (\ref{frmSO}) that $DF$ and $DG$ are of the form
\begin{equation*}
\fl
DF=
 \left(
  \begin{array} {cc}
    {2 d \xi^* + f_{11} \eta^* \lambda_1^k  + 3 s_{03}\lambda_1^k (\xi^*)^2}  & {c  + f_{11} \xi^* \lambda_1^k} \\
    {b \lambda_1^{-k}\lambda_2^m + 2 \lambda_2^m l_{02} \xi^*} & a  \lambda_2^m
  \end{array}
  \right),
\end{equation*}
\begin{equation*}
\fl
DG=
 \left(
  \begin{array}{cc}
    a \lambda_2^m  &  {b \lambda_1^{-k} \lambda_2^m + 2 \lambda_2^m l_{02} \eta^*}\\
    {c  + f_{11} \eta^* \lambda_1^k} & {2 d \eta^* + f_{11} \xi^* \lambda_1^k  + 3 \lambda_1^k s_{03} (\eta^*)^2}
  \end{array}
 \right)
\end{equation*}
plus terms of order $O(k\lambda_1^{2k})$.
Now we compute the Jacobian as
\begin{equation}
\label{jkm}
\begin{array}{l}
{\ds J(T_{km})\bigl|_{(\xi^*,\eta^*)} = \det(DF^{-1}DG)=\frac{\det(DG)}{\det(DF)}= } \\[1.4ex]
{\ds \frac{-bc\lambda_1^{-k}\lambda_2^m + 2 ad \lambda_2^m \xi^* -
b f_{11}\lambda_2^m \xi^* - 2 c l_{02} \lambda_2^m
\xi^*+o(\lambda_1^k)}{-bc\lambda_1^{-k}\lambda_2^m + 2 ad
\lambda_1^k\lambda_2^m \eta^* - b f_{11}\lambda_2^m \eta^* - 2 c
l_{02} \lambda_2^m \eta^*+o(\lambda_1^k)}.}
\end{array}
\end{equation}
When the relation (\ref{keqm}) is fulfilled we can rewrite (\ref{jkm}) as
\[
J(T_{km})\bigl|_{(\xi^*,\eta^*)} = \frac{-bc + Q\lambda_1^k\;\xi^*
+o(\lambda_1^k)}{-bc + Q\lambda_1^k\;\eta^* +o(\lambda_1^k)}
\]
that gives relation (\ref{FRMJac}).

\subsection{Proof of Lemma~\ref{lmgenel} }
\label{se:birkhoff}
Due to the reversibility, we can prove Lemma~\ref{lmgenel}
directly for the truncated map ${\ds H: \ \bar x = \Mt + \ct x -
y^2, \ \bar y = -\frac{1}{\ct} \Mt + \frac{1}{\ct} y +
\frac{1}{\ct} \bar x^2}$. We will use the following facts for this
map: it can be written in the explicit form (\ref{frmex}) and that
for $M>-\frac{1}{4}(\ct -1 )^2$ it has a pair of symmetric fixed
points $P^{+}=(p_+,p_+)$ and $P^{-}=(p_{-},p_{-})$ for which
coordinates the formula (\ref{p+p-}) holds.
Denote by $p$ either $p_+$ or $p_{-}$ and let us assume that the
corresponding fixed point $P$ (i.e. $P^{+}$ or $P^{-}$) is
elliptic. Then, $\ct$ and $\Mt$ have to take values from the open
regions in the $(\ct,\Mt)$-space of parameters given in
Figure~\ref{Figdiag31}.

The first step in our process is to shift
the new origin of coordinates into the point $(p,p)$ and to
perform (Jordan) linear normal form, which leads our map to the
following form
\begin{equation}
\label{eq:maprot}
\begin{array}{l}\!\!\!\!\!\!\!\!\!\!\!\!\!\!\!\!\!\!\!\!\!\!\!\!\!\!\!\!\!\!\!\!\!\!\!\!\!\!\!\!
\displaystyle \bar{x} =  \cos\psi \cdot x - \sin\psi\cdot y -
\frac{2p\cos\psi}{\ct \sin\psi}y^2 +  \\
\!\!\!\!\!\!\!\!\!\!\!\!\!\!\!\displaystyle +
  \frac{1-4p^2- \ct \cos\psi}{4 \ct^2 p^2 \sin\psi}\left(-\ct \sin\psi\cdot x + (1-\ct \cos\psi) y + 2p y^2 \right)^2
   , \\ \\
\!\!\!\!\!\!\!\!\!\!\!\!\!\!\!\!\!\!\!\!\!\!\!\!\!\!\!\!\!\!\!\!\!\!\!\!\!\!\!\!
\displaystyle \bar{y} = \sin\psi \cdot x + \cos\psi\cdot y -
\frac{2p}{\ct} y^2
+ \\
\!\!\!\!\!\!\!\!\!\!\!\!\!\!\!\displaystyle  +  \frac{1}{4 \ct
p^2}\left(-\ct \sin\psi\cdot x + (1-\ct \cos\psi) y + 2p y^2
\right)^2.
\end{array}
\end{equation}
where $x,y,\bar{x},\bar{y}$ stand again for the new variables. The
linear part of (\ref{eq:maprot}) is a rotation of angle $\psi$.
At that point, we will lead our map into the so-called Birkhoff
Normal Form up to order $3$:  $\;\bar{z}= \rme^{\rmi \psi} z +
d_{21} z^2 z^* + \mathcal{O}_4$. To do it,
we need to assume that $\lambda=\rme^{\rmi \psi}$ is not a
$k$th-root of unity for $k=3,4$ (the cases $k=1,2$ correspond to
parabolic fixed points and, thus, respond for boundaries of
existence regions of elliptic fixed points). The coefficient
$B_1\equiv -i d_{21} e^{-i\psi}$ is called the first Birkhoff
coefficient. By the Arnol'd-Moser Twist Theorem~\cite{SM95},
the inequality $B_1\neq 0$ (together with the
absence of strong resonances) ensures that the elliptic point is
generic or, in other words, KAM-stable.

Introducing complex coordinates $z=x+\rmi y$, $z^*=x-\rmi y$, map
(\ref{eq:maprot}) takes the form
\[
\bar z = {\rme}^{\rmi {\psi}} z + A_{20}z^2 + A_{11}zz^* +
A_{02}(z^*)^2 + A_{21}z^2z^* + \mathcal{O}_4(z,z^*).
\]
Since we are assuming $\rme^{\rmi \psi}$ not to be a $3$rd or
$4$th root of unity (and also $\psi\neq 0,\pi$), our map can be
lead into BNF up to order $3$ and, afterwards, provides the
following formula for the first Birkhoff coefficient
\begin{equation}
\label{eq:B1_H} B_1=\frac{(\ct+1-2p)(\ct+1+2p)(\ct-1+2p)}{32 \ct^4
p \sin^3 \psi (2\cos \psi +1)} \  P_4(\ct,p),
\end{equation}
where
\begin{eqnarray*}
P_4(\ct,p) & = 64p^4 + 8(1-\ct)p^3 - 4 (3\ct^2+4\ct+3) p^2 + \\
           & 2(\ct-1)(\ct+1)^2 p - (\ct-1)^2 (\ct+1)^2.
\end{eqnarray*}

\begin{figure}[htb]
\begin{center}
\includegraphics[width=13cm]{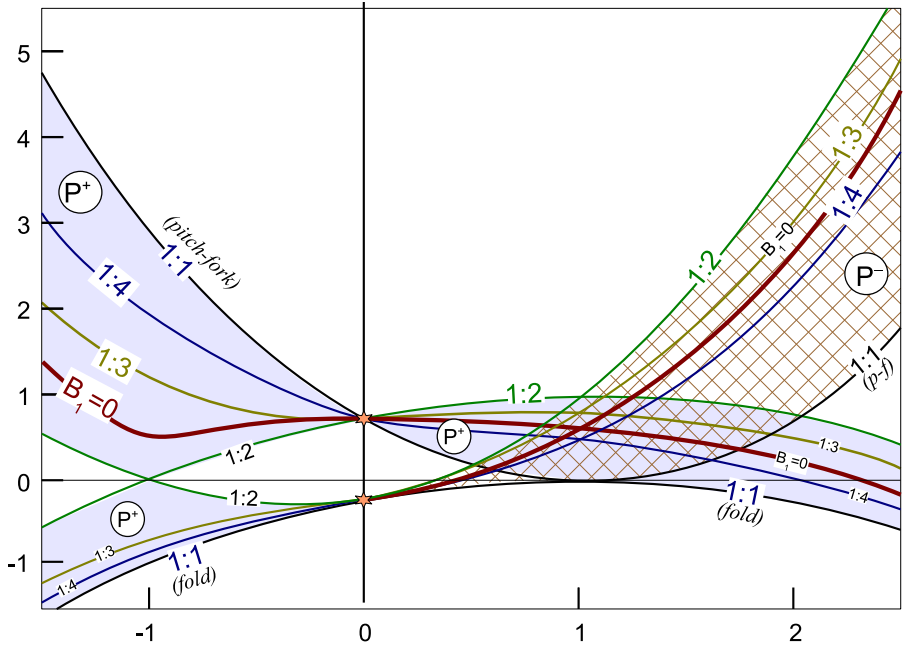} 
\caption{In the $(\ct,\Mt)$-plane, three grey and
one hatching regions correspond to the existence of elliptic
points: $P^+$ for the grey regions and $P^-$ for the hatching one.
Lines corresponding to the main resonances and vanishing the first
Birkhoff coefficient for the elliptic point are shown and
labelled.}
\label{fig:b1=0}
\end{center}
\end{figure}

Using this formula\footnote{Notice that in the particular case
$\ct=-1$, map $H$ corresponds to $\mathcal{H}^2$, where $\;\;
\mathcal{H}: \ \bar{x}=y, \; \bar{y} = M - x - y^2\;\;$  is the
H\'enon map. In this case we have  that $B_1^{H}(\psi) =
2B_1^{\mathcal{H}(\varphi)}$ where $\psi=2\varphi$. It is not
hard to check now that, for a fixed point of $\mathcal{H}$ with
$p=-\cos\varphi$,  the following relation holds
\[
B_1^{\mathcal{H}^2}= B_1^{\mathcal{H}^2}=
\frac{1}{4\sin^2\varphi}\cdot\frac{(1+\cos\varphi)(1+4\cos\varphi)}{\sin\varphi(1+2\cos\varphi)},
\]
which differs from the well-known formula for the H\'enon map (see
e.g. \cite{B87}) only by the non-zero factor $\frac14
\sin^{-2}\varphi$.} we represent in Figure~\ref{fig:b1=0} curves
$B_1(\ct,\Mt)=0$, where the elliptic fixed point can be, a priori,
not KAM-stable. In this Figure curves related to the strong
resonances are presented. For resonances 1:1 and 1:2 the
equations of the corresponding curves are given in
Section~\ref{sec4}. The equations of the 1:3 and 1:4 curves
are as follows:
$$
\Mt = \frac{\ct^2+1}{4} \pm \frac{\sqrt{\ct^2+1}}{2}(1-\ct)\;\;
\mbox{for 1:4 resonance}
$$
and
$$
\Mt = \frac{\ct^2+\ct+1}{4} \pm \frac{\sqrt{\ct^2+\ct+
1}}{2}(1-\ct)\;\; \mbox{for 1:3 resonance}.
$$
This completes the proof.

\vskip0.25cm

{\bf Acknowledgements.} The authors thank D.Turaev and L.Shilnikov
for fruitful discussions and remarks. The second, third and fifth authors have been supported partially by the RFBR grants No.07-01-00566, No.07-01-00715 and
No.09-01-97016-r-povolj'e. The first and fourth authors have been partially supported by the MICIIN/FEDER grant number MTM2009-06973
and by the Generalitat de Catalunya grant number 2009SGR859.
S.~Gonchenko thanks Centre de Recerca Matem\`atica (Bellaterra) for very nice hospitality in 2008.

\end{document}